\documentclass{elsart}
\usepackage[english]{babel}
\usepackage{latexsym}
\usepackage{fancyhdr}
\usepackage{stmaryrd}
\usepackage{amsmath,amssymb,graphicx}
\usepackage[numbers]{natbib}
\journal{Stochastic Processes and their Applications,}
\def\cite{\citet}

\def \R{\mathbb{R}}
\def \D{\mathbb{D}}

\newcommand{\E}{\mathbb{E}}
\newcommand \Ek{\E_{t_k}}

\begin{document}
\begin{frontmatter}
  
  \title{Error expansion for the discretization of Backward Stochastic Differential Equations}
  \author[Grenoble]{Emmanuel Gobet},
  \ead{emmanuel.gobet@imag.fr}
  \author[Paris]{C{\'e}line Labart}
  \ead{labart@cmapx.polytechnique.fr}
  \address[Grenoble]{ ENSIMAG - INPG, IMAG - LMC, BP 53, 38041 Grenoble cedex 9, FRANCE}
  \address[Paris]{CMAP,Ecole Polytechnique, Route de Saclay, 91128 Palaiseau, FRANCE}
  
  \begin{abstract}
    We study the error induced by the time discretization of a decoupled forward-backward stochastic differential equations $(X,Y,Z)$. The forward component $X$ is the solution of a Brownian stochastic differential equation and is approximated by a Euler scheme $X^N$ with $N$ time steps. The backward component is approximated by a backward scheme. Firstly, we prove that the errors $(Y^N-Y,Z^N-Z)$ measured in the strong $L_p$-sense ($p \geq 1$) are of order $N^{-1/2}$ (this generalizes the results by \cite{zhang_04}). Secondly, an error expansion is derived: surprisingly, the first term is proportional to $X^N-X$ while residual terms are of order $N^{-1}$.
  \end{abstract}

\end{frontmatter}

\section{Introduction}
Let $(\Omega, \mathcal{F},\mathbb{P})$ be a given probability space on which is defined a $q$-dimensional standard Brownian motion $W$, whose natural filtration, augmented with $\mathbb{P}$-null sets, is denoted by $(\mathcal{F}_t)_{0\le t\le T}$ ($T$ is a fixed terminal time). We consider the solution $(X,Y,Z)$ to a decoupled forward-backward stochastic differential equation (FBSDE in short). Namely, $X$ is the $\mathbb{R}^d$-valued process solution of 
\begin{equation} \label{eqs4}
  X_t= x + \int_0^t b(s,X_s) ds + \int_0^t \sigma(s,X_s)dW_s, 
\end{equation}
and $Y$ (resp. $Z$) is a real-valued adapted (resp. predictable $\mathbb{R}^q$-valued) process solution of 
\begin{equation} \label{eqs1}
  -dY_t=f(t,X_t,Y_t,Z_t)dt-Z_tdW_t, \;\; Y_T=\Phi(X_T).
\end{equation}
We assume standard Lipschitz properties on the coefficients, which ensure existence and uniqueness in appropriate $L_2$-spaces (see \cite{pard:peng:90}, or \cite{ma:yong:99} for numerous references). During the last decade, more and more attention has been paid to these equations, because of their natural applications in Mathematical Finance or in the probabilistic resolution of semi-linear partial differential equations (PDE in short): see \cite{karoui_97} or \cite{pard:98}.

Our aim is to study the most usual time approximation of $(X,Y,Z)$. For $X$, we use the Euler scheme $X^N$ with $N$ discretization times $(t_k=kh)_{0 \leq k \leq N}$ ($h=\frac{T}{N}$ is the time step). For convenience, set $\Delta W_k= W_{t_{k+1}}-W_{t_k}$ ($\Delta W^l_k$ component-wise). $X^N$ is defined by $X^N_0=x$ and
\begin{equation}
  \label{eq:euler}
t\in[t_k,t_{k+1}],\  X^N_{t}=X^N_{t_k}+b(t_k,X^N_{t_k}) (t-t_k) + \sigma(t_k,X^N_{t_k}) (W_t- W_{t_k}).
\end{equation}
The backward SDE \eqref{eqs1} is approximated by $(Y^N,Z^N)$ defined in a backward manner by $ Y_{t_N}^N= \Phi(X^N_{t_N})$ and
\begin{align}\label{eqs2}
  Y^N_{t_k}  = & \mathbb{E}_{t_k}(Y^N_{t_{k+1}}) +
  h\mathbb{E}_{t_k}f(t_k,X^N_{t_k},Y^N_{t_{k+1}},Z^N_{t_k}),\\
  \label{eqs3} hZ^N_{t_k} =& \mathbb{E}_{t_k}(Y^N_{t_{k+1}} \Delta W_k^*),
\end{align}
where $\mathbb{E}_{t_k}$ is the conditional expectation w.r.t. $\mathcal{F}_{t_k}$ and $*$ is the transpose operator. Additional tools are needed to derive a fully implementable scheme, in particular for the computations of conditional expectations. We refer to \cite{bouc:touz:04} for Malliavin calculus techniques, or to \cite{lemor_05,lemor_06} for empirical regression methods. In this work, we leave these further questions and we only address the error analysis between $(Y,Z)$ and $(Y^N,Z^N)$. \\
On the one hand, \cite{zhang_04} proves (in a slightly different form) that the error $\max_{k\le N}\|Y^N_{t_k}-Y_{t_k}\|_{L_2}\leq C N^{-1/2}$. This is done under rather minimal Lipschitz assumptions on $b,\sigma,f,\Phi$. On the other hand, when $f$ does not depend on $z$ and the coefficients are smooth, one knows that $|Y^N_{0}-Y_{0}|\le CN^{-1}$ (see \cite{chevance_97}). We aim at filling the gap regarding these two different rates of convergence. In the following, we prove that
\begin{itemize}
\item the Chevance's results extend to the case of $f$ depending also on $z$.
\item the rate $N^{-1}$ holds true also for the difference $|Z^N_{0}-Z_{0}|$.
\item more generally, for the other discretization times $t_k$, we expand the error as $$\big|Y^N_{t_k}-Y_{t_k}-\alpha_k \cdot (X^N_{t_k}-X_{t_k})\big|\le C N^{-1}\lor |X^N_{t_k}-X_{t_k}|^2$$ (for an explicit and bounded random vector $\alpha_k$). 
\item an analogous expansion is available for $Z$.
\end{itemize} 
Since $|X^N_{t_k}-X_{t_k}|^2$ has the same order in $L_p$ than $N^{-1}$, the error on $Y$ is mainly due to the error $X^N_{t_k}-X_{t_k}$. Thus, Zhang's results are a consequence of this expansion, and Chevance's ones as well since $X^N_0=X_0$. The gap is filled.\\
In addition, we learn from this expansion that if one could perfectly simulate $X$
(as for Brownian motion with constant drift, geometric Brownian motion or
Ornstein-Uhlenbeck process), the error on the BSDE would be of order $N^{-1}$ and not
$N^{-1/2}$ as stated by Zhang's results. Also, if one could use a discretization
scheme for $X$ of order 1 for the strong error (for instance Milshtein scheme
whenever possible), the error on the BSDE would be of order $N^{-1}$ (we would need
to extend our analysis to other discretization schemes, this is straightforward for
the Milshtein scheme).

The paper is organized as follows. In Section \ref{section1}, we define the
assumptions on the coefficients, recall the connection between BSDEs and semi-linear
PDEs (which is important for our analysis). Finally, we state our main
results. Firstly in Theorem \ref{theo1}, we extend Zhang's results to $L_p$
norm. Secondly in Theorem \ref{theo2}, we expand the error on $Y$. Lastly in
Theorem \ref{theo4}, we deal with the error on $Z$. Naturally, stronger and stronger
assumptions are required for theses theorems. Proofs of the three results are
postponed to Sections \ref{section:proof:theo1}, \ref{section:proof:theo2} and
\ref{section:proof:theo4}: we combine BSDE techniques, martingale estimates and Malliavin calculus.

{\bf Notation.}
\begin{itemize}
\item {\bf Differentiation.} If $g: \mathbb{R}^d \mapsto \mathbb{R}^q$ is a differentiable function, its
  gradient $\nabla_x g(x)=(\partial_{x_1} g(x),...,\partial_{x_d}g(x))$ takes
  values in $\mathbb{R}^q \otimes \mathbb{R}^d$. At many places, $\nabla_x g(x)$ will
  simply be denoted $g'(x)$.  If $g: \mathbb{R}^d \mapsto \mathbb{R}$ is a twice differentiable function,
  its Hessian $H_x(g)$  takes
  values in $\mathbb{R}^d \otimes \mathbb{R}^d$: $(H_x(g))_{i,j}=\partial^2_{x_i x_j}  g$.
  If $g:\mathbb{R}^d \times \mathbb{R}^q \mapsto \mathbb{R}$, $g^{''}_{xy}(x,y)$
  takes values in $\mathbb{R}^{d}\otimes\R^q$: $(g^{''}_{xy})_{ij}=\frac{\partial^2
    g}{\partial x_i \partial y_j},$ for $1 \leq i \leq d, 1 \leq j \leq q$.
\item{\bf Function spaces.} For an integer $k \geq 1$, we denote by $C_b^{k/2,k,k,k}$ the set of continuously differentiable functions $\phi:(t,x,y,z)\in [0,T]\times \R^d\times\R\times\R^q\mapsto \phi(t,x,y,z)$ such that the partial derivatives $\partial_t^{l_0}\partial^{l_1}_{x}\partial_y^{l_2}\partial_z^{l_3}\phi(t,x,y,z)$ exist for $2l_0+l_1+l_2+l_2\le k$ and are uniformly bounded. The analogous set of functions that not depend on $y$ and $z$ is denoted by $C_b^{k/2,k}$. This set is denoted by $C_b^{(k+\alpha)/2,k+\alpha}$ ($\alpha\in]0,1[$) if in addition the highest derivatives are H\"older continuous with index $\alpha$ w.r.t. $x$ and $\alpha/2$ w.r.t. $t$ (for a precise definition, see 
\cite{lady:solo:ural:68}).
\item{\bf Norm.} For a $d$-dimensional vector $U$, we set $|U|^2=\sum_{i=1}^d
  U_i^2$. For a $d\times q$-dimensional matrix $A$, $A_i$ denotes its {\it i}-th column, and $A^i$ its {\it i}-th row. Moreover, $|A|^2=\sum_{i,j=1}^{d,q} A_{i,j}^2$. 
\item {\bf Constants.} Let $C$ denote a generic constant which may depend on the coefficients $b,\sigma,f,\Phi$ and on the dimensions $d$ and $q$. We will keep the same notation $K(T)$ for all finite, nonnegative, and nondecreasing functions w.r.t. $T$: they do not depend on $x$ and $h$. The generic notation $K(T,x)$ stands for any function
  bounded by $K(T)(1+|x|^q)$, for some $q \geq 0$.
\item {\bf \mathversion{bold}$O(U)$ and $O_k(h)$.} A random vector $R$ is such that $R=O(U)$ for a nonnegative random variable $U$ if $|R| \leq K(T,x)U$ (in particular, $R=O(h)$ means $|R| \leq K(T,x)h$). The notation $R=O_k(h^p)$ means $|R| \leq \lambda^N_k
   h^p$, where $\lambda^N_k$ is $\mathcal{F}_{t_k}$-measurable,
   $\sup_N\E(\sup_k|\lambda^N_k|^q) \le K(T,x)$,  for $q\ge 1$.
\item{\bf \mathversion{bold}$\mathbb{E}_{t_k}$ and $\mbox{Var}_{t_k}$.} $\mathbb{E}_{t_k}$ is the
  conditional expectation w.r.t. $\mathcal{F}_{t_k}$ and
  $\mbox{Var}_{t_k}(X)=\E_{t_k}(X^2)-(\E_{t_k}(X))^2$.
\item{\bf Malliavin calculus}. We use the notations of \cite{nualart_95} for weak spaces $\mathbb{D}^{k,p}$.
\item{\bf Discretization} Let $s \in [t_k, t_{k+1}[$. We define $\eta(s)=t_k$.
\end{itemize}

\section{Main results}
\label{section1}
\subsection{Hypotheses}
The coefficients $b$ : $[0,T] \times \mathbb{R}^d
\rightarrow \mathbb{R}^d$, $\sigma$ : $[0,T] \times \mathbb{R}^d
\rightarrow \mathbb{R}^{d\times q}$, $f : [0,T] \times \mathbb{R}^d  \times \mathbb{R} \times
\mathbb{R}^q \rightarrow \mathbb{R}$ and $\Phi :\mathbb{R}^d \rightarrow
\mathbb{R}$ satisfy one of the following set of assumptions.

\begin{hypo}\label{hyp3}
  The functions $b, \sigma, f$ and $\Phi$ are bounded in $x$, are uniformly Lipschitz continuous w.r.t. $(x,y,z)$ and H{\"o}lder continuous of parameter $\frac{1}{2}$ w.r.t. $t$. In addition, $\Phi$ is of class $C_b^{2+\alpha}$ for some $\alpha\in]0,1[$ and the matrix-valued function $a=\sigma \sigma^*$ is uniformly elliptic.
\end{hypo}

\begin{hypo}\label{hyp4}
  Assume Hypothesis \ref{hyp3} and that the functions $b$, $\sigma$ are in
    $C^{\frac{3}{2},3}_b$, $f$ is in
    $C^{\frac{3}{2},3,3,3}_b$, $\Phi$ is in $C_b^{3+\alpha}$ for some $\alpha\in]0,1[$.
\end{hypo}

\begin{hypo}\label{hyp5}
  Assume Hypothesis \ref{hyp3} and that the functions $b$, $\sigma$ are in
    $C^{2,4}_b$, $f$ is in
    $C^{2,4,4,4}_b$, $\Phi$ is in $C^{4+ \alpha}$ for some $\alpha\in]0,1[$.
\end{hypo} 
We do not assert that these smoothness and boundedness conditions are the weakest ones for our error analysis, but they are sufficient. Investigations regarding minimal assumptions would be certainly interesting but it is beyond the scope of the paper.

\subsection{Connection between Markovian BSDE's and semi-linear parabolic PDE's}
We recall classical results connecting $(Y,Z)$ and the solution and its gradient of the following semi-linear PDE on $[0,T]\times \mathbb{R}^d$:
\begin{align}\label{EDP}
 &(\partial _t + \mathcal{L}_{(t,x)})u(t,x)+f \left(t,x,u(t,x), \nabla_x u(t,x) \sigma(t,x)\right)=0,\\
&u(T,x)=\Phi(x),\notag
\end{align}
where $\mathcal{L}_{(t,x)}$ is the second order differential operator
$$
\mathcal{L}_{(t,x)}=\frac{1}{2}\sum_{i,j} [\sigma \sigma^*]_{ij}(t,x)\partial_{x_ix_j}^2 +\sum_i
b_i(t,x)\partial_{x_i}$$ 
(see for instance \cite{ma_02} or \cite{pard:98}).
\begin{prop}\label{theo5}Under Hypothesis \ref{hyp3}, one has 
\begin{align}\label{eqs41}
  \forall t\in [0,T],\;
  Y_t=u(t,X_t), \; Z_t= \nabla_x u(t,X_t) \sigma(t,X_t),
\end{align}
where $u$ is the unique classic solution $C^{1,2}_b$ of the PDE \eqref{EDP}.

In addition under Hypothesis \ref{hyp4}, $u \in C^{\frac{3}{2},3}_b$, and under Hypothesis \ref{hyp5}, $u \in C^{2,4}_b$.
\end{prop}
The first result of this Proposition corresponds to Theorem $2.1$ of \cite{delarue_05}. The two last regularity results can be proved in the same way. In fact for this, we would only need $b,\sigma$ to be in $C^{1+\alpha/2,2+\alpha}_b$; the additional smoothness is used later for Malliavin calculus computations.

\subsection{Main results}
We now turn to the statement of our results. Remind the following well-known upper bound on the Euler Scheme, which is useful in the sequel.
\begin{prop}\label{theo3}
  Let $\sigma$ and $b$ be Lipschitz continuous. Then
  \begin{equation*}
    \forall p \geq 1, [\E(\sup_{t \leq T}|X^N_t-X_t|^p)]^{\frac{1}{p}} \leq K(T,x)\frac{1}{\sqrt{N}}.
  \end{equation*}
\end{prop}
In fact, for all $p\ge 1$ one has 
  \begin{equation}\label{remarque1}
     [\E_{t_i}(\sup_{t_i \leq t \leq T}|X^N_t-X_t|^p)]^{\frac{1}{p}}
    \leq K(T,X_{t_i})\frac{1}{\sqrt{N}}+| X^N_{t_i}-X_{t_i}|.
  \end{equation}
Our first result is an extension of the $L_2$ estimates in \cite{zhang_04} to $L_q$ estimates (see also \cite{lemor_05}).
\begin{thm}\label{theo1}
  Let us assume Hypothesis \ref{hyp3}.
  Let $q>0$. We define the error 
  $$
  e_q(N)=\big[\max_{0 \leq k \leq N} \mathbb{E}|Y_{t_k}-Y_{t_k}^N|^q + \mathbb{E}
  (\sum_{k=0}^{N-1} \int_{t_k}^{t_{k+1}} |Z^N_{t_k}-Z_t|^2 dt)^{\frac{q}{2}}\big]^{\frac{1}{q}},
  $$
  where $Y^N$ and $Z^N$ are defined by (\ref{eqs2}) and (\ref{eqs3}). Then
  $|e_q(N)|\leq K(T,x) \frac{1}{\sqrt{N}}$.
\end{thm}

By slightly strengthening the smoothness assumptions on $b,\sigma, f$ and $\Phi$, we are able to expand the error on $Y$.
\begin{thm}\label{theo2}
  Let us assume Hypothesis \ref{hyp4}. Then, the following expansion holds
  \begin{align*}
    Y^N_{t_k} - Y_{t_k} = & \nabla_x u(t_k,X_{t_k})( X^N_{t_k}-X_{t_k}) +O_k(\frac{1}{N}) + O(
    |X^N_{t_k}-X_{t_k}|^2).
  \end{align*}
\end{thm}
In view of Proposition \ref{theo3}, $|X^N_{t_k}-X_{t_k}|^2$ and $N^{-1}$ have the same
order (in $L_p$). Hence it turns out that $ \nabla_x u(t_k,X_{t_k})( X^N_{t_k}-X_{t_k})$ is the first order term in the error $Y^N_{t_k} - Y_{t_k}$. Obviously, this estimate implies that of Theorem \ref{theo1}. As mentioned in the introduction, the evaluation of $Y_0$ by $Y^N_0$ has still an accuracy of order $N^{-1}$ since initial conditions for $X^N$ and $X$ coincide. Note that if there is no discretization error for the process $X$, $Y^N_{t_k} - Y_{t_k}=O(\frac{1}{N})$, a fact which is not clear from equations \eqref{eqs2} and \eqref{eqs3}. A nice situation corresponds to $\sigma$ independent of $x$ (this is a very specific situation where Euler and Milshtein schemes are equal): in that case $\|X^N_{t_k} - X_{t_k}\|_{L_p}=O(N^{-1})$ and one gets the order of accuracy $N^{-1}$ for $Y$.

For $Z$ which plays the role of a gradient relatively to $Y$, we get an analogous result about the error, up to increasing by 1 the degree of smoothness of the coefficients.
\begin{thm}\label{theo4}
  Let us assume Hypothesis \ref{hyp5}. Then, the following expansion holds
  \begin{align*}
    Z^N_{t_k}-Z_{t_k}=  \big(\nabla_x[\nabla_x u \ \sigma]^*(t_k, X_{t_k})
    (X^N_{t_k}-X_{t_k})\big)^* + O_k(\frac{1}{N}) + O(|X^N_{t_k}-X_{t_k}|^2).
  \end{align*}
\end{thm}
\begin{rem}\label{remarque5} The above results are sufficient to derive the weak convergence
  of the renormalized error process $[\sqrt{N}(Y^N_t-Y_t)]_{0 \leq t \leq
    T}$ and $[\sqrt{N}(Z^N_t-Z_t)]_{0 \leq t \leq T}$, except that one has to define $Y^N$ and $Z^N$ between discretization times. For $t \in [t_k, t_{k+1}[,$ analogously to (\ref{eqs2}) and (\ref{eqs3}) we define
  \begin{align*}
    Y^N_t&=\E_t\big(Y^N_{t_{k+1}}+(t_{k+1}-t)f(t,X^N_t,Y^N_{t_{k+1}},Z^N_t)\big),\\
    Z^N_t&=\frac{1}{t_{k+1}-t}\E_t\big(Y^N_{t_{k+1}}(W_{t_{k+1}}-W_t)^*\big).
  \end{align*}
  Theorems \ref{theo2} and \ref{theo4} can be extended to all $t \in [0,T]$. We have
   \begin{align*}
    Y^N_t - Y_t = & \nabla_x u(t,X_t)( X^N_t-X_t) +O_t(\frac{1}{N}) + O(
    |X^N_t-X_t|^2),\\
    Z^N_t-Z_t=&  \big(\nabla_x[\nabla_x u \ \sigma]^*(t,X_t)
    (X^N_t-X_t)\big)^* + O_t(\frac{1}{N}) + O(|X^N_t-X_t|^2).
  \end{align*}
  Theorem $3.5$ of \cite{kurtz_91} allows us to establish the weak convergence
  of the processes $\sqrt{N}(Y^N-Y)$, and $\sqrt{N}(Z^N-Z)$. Indeed, the process $[\sqrt{N}(X^N_t-X_t)]_{0 \leq t \leq T}$ weakly converges to the solution of
  \begin{align*}
    U_t=&\sum_{i=1}^q\int_0^t
    \nabla_x \sigma_i(s,X_s)U_s dW^i_s+ \int _0^t \nabla_x b(s,X_s)U_s ds\\
    &+\frac{1}{\sqrt{2}}\sum_{i,j=1}^q\int_0^t \sum_{k=1}^d \partial_{x_k} \sigma_i(s,X_s)\sigma_{kj}(s,X_s)dV_s^{ij},
  \end{align*}
   where $(V^{ij})_{1\leq i,j \leq q}$
  are independent standard Brownian motions and independent of
  $W$. Furthermore, the convergence is stable (see \cite{jacod_98}). Hence, $[\sqrt{N}(X_t^N-X_t),\sqrt{N}(Y_t^N-Y_t),\sqrt{N}(Z_t^N-Z_t),X_t]_{0\le t\le T}$ weakly converges to $[(U_t, \nabla_x u(t,X_t)U_t,( [\nabla_x[\nabla_x u \ \sigma]^*(t, X_t)] U_t)^*, X_t]_{0 \leq t \leq T}$.
  \end{rem}
  \vspace{-0.5cm}
\section{Proof of theorem \ref{theo1}}\label{section:proof:theo1}
\vspace{-0.5cm}{\bf Extra notations for all the proofs.} For any process $U$ (except the Brownian increments $\Delta W_k$), we define $\Delta U_k=
U^N_{t_k}-U_{t_k}$ . 
Let $\theta_s$ denote $(s,X_s,Y_s,Z_s)$ and $f^N_{t_k}$ denote $f(t_k,X^N_{t_k},Y^N_{t_{k+1}},Z^N_{t_k})$. \\
$\overline{Z}_{t_k}$ is defined as 
$h\overline{Z}_{t_k} := \mathbb{E}_{t_k} \int_{t_k}^{t_{k+1}} Z_s ds$ and we put $\Delta\overline{Z}_k=Z_{t_k}^N-\overline{Z}_{t_k}$.

If $q=2$, the result has already been proved in \cite{lemor_05}, under Lipschitz conditions on $b, \sigma, f, \Phi$. Thanks to the inequality $\mathbb{E} |U|^q \leq (\mathbb{E}
|U|^{2p})^{\frac{q}{2p}}$ for $2p\geq q$, we only need to prove the theorem for
$q=2p$, where $p \in \mathbb{N}^*$.\\
First, we give some estimates which can be easily established. We have, under Hypothesis
\ref{hyp3}, $\forall s \in [t_k, t_{k+1}]$,
\begin{align}\label{eqs18}
  \mathbb{E}_{t_k} ( |X_s-X_{t_k}|^{2p}+|Y_s-Y_{t_k}|^{2p} + |Z_s-\overline{Z}_{t_k}|^{2p} )\leq Ch^p.
\end{align}
In the following computations, these estimates are repeatedly used. 
\subsection{Proof of  $ \max_{0 \leq k \leq N}
  \mathbb{E}|Y_{t_k}-Y_{t_k}^N|^{2p}=O(h^{p})$.}\vspace{-0.5cm}
We prove the following result, which is a bit more general.
\begin{prop}\label{prop2}
  $\max_{i \leq k \leq N}
  \E_{t_i}|Y_{t_k}-Y_{t_k}^N|^{2p}=O_i(h^p)+|\Delta X_i|^{2p}$.
\end{prop}
By taking $i=0$, we get $ \max_{0 \leq k \leq N} \mathbb{E}|Y_{t_k}-Y_{t_k}^N|^{2p}=O(h^{p})$.\\
{\bf Assume that we have}
\begin{equation}
  \label{eq:y}
  |\Delta Y_k|^2 \leq \displaystyle (1+ Ch)\mathbb{E}_{t_k}|\Delta
  Y_{k+1}|^2+Ch|\Delta X_k|^2 +Ch^2.
\end{equation}Then, using the inequality $(a+b)^p \leq
a^p(1+ \epsilon(2^{p-1}-1)) + b^p(1+\frac{2^{p-1}-1}{\epsilon^{p-1}})$ for $0<\epsilon<1$, we deduce
\begin{equation*}
  |\Delta Y_k|^{2p} \leq \displaystyle (1+ Ch)^{p+1}\mathbb{E}_{t_k}|\Delta
  Y_{k+1}|^{2p} + C^ph^p(|\Delta X_k|^2 +Ch)^p(1+\frac{C}{h^{p-1}}).
\end{equation*}
Take the conditional expectation w.r.t. $\mathcal{F}_{t_i}$ to get
$
\mathbb{E}_{t_i}|\Delta Y_k|^{2p} \leq (1+C h) \mathbb{E}_{t_i}|\Delta
Y_{k+1}|^{2p}+h(h^p+\E_{t_i}|\Delta X_k|^{2p}).
$
Using \eqref{remarque1} for $|\Delta X_k|$ and Gronwall's lemma yields $\max_{i \leq k \leq N}
\mathbb{E}_{t_i}|Y_{t_k}-Y_{t_k}^N|^{2p}=O_i(h^p) + |\Delta X_i|^{2p}.$\qed\\
{\bf Now we prove the inequality \eqref{eq:y}.} From (\ref{eqs1}) and (\ref{eqs2}) we obtain
\vspace{-0.3cm}
\begin{equation}\label{eqs0}
  \Delta Y_k = \mathbb{E}_{t_k}(\Delta Y_{k+1}) + \mathbb{E}_{t_k}
  \int_{t_k}^{t_{k+1}}(f^N_{t_k} -f(\theta_s)) ds.
\end{equation}
By applying Young's inequality, that is $(a+b)^2 \leq (1+ \gamma h)a^2 +
(1+\frac{1}{\gamma h})b^2$, where $\gamma$ will be fixed later, and using the Lipschitz property of $f$, we get 
\begin{align}\label{eqs5}
  |\Delta Y_k|^2& \leq (1+ \gamma h)(\mathbb{E}_{t_k}(\Delta Y_{k+1}))^2+
  C( h+ \frac{1}{\gamma }) [ h^2+
  \mathbb{E}_{t_k}\int_{t_k}^{t_{k+1}}
  |X_s-X^N_{t_k}|^2 ds]\notag\\
  & + C( h+ \frac{1}{\gamma })
  [\mathbb{E}_{t_k}\int_{t_k}^{t_{k+1}} |Y_s-Y^N_{t_{k+1}}|^2 ds +\mathbb{E}_{t_k} \int_{t_k}^{t_{k+1}} |Z_s-Z^N_{t_k}|^2 ds].
\end{align}
Let us introduce $\overline{Z}_{t_k}$ (see extra notations at the beginning of Section \ref{section:proof:theo1}):
\begin{equation}\label{eqs6}
  \mathbb{E}_{t_k}   \int_{t_k}^{t_{k+1}} |Z_s-Z^N_{t_k}|^2 ds = \mathbb{E}_{t_k}
  \int_{t_k}^{t_{k+1}} |Z_s-\overline{Z}_{t_k}|^2 ds + h \mathbb{E}_{t_k}|\overline{Z}_{t_k}-Z_{t_k}^N|^2. 
\end{equation}
Thanks to the Cauchy Schwarz inequality we have 
$$
|\mathbb{E}_{t_k}(\Delta Y_{k+1}\Delta W^l_k)|^2 \leq h\{
\mathbb{E}_{t_k}(|\Delta Y_{k+1}|^2) -
|\mathbb{E}_{t_k}(\Delta Y_{k+1})|^2\}.
$$
Hence, as $h\overline{Z}_{t_k} =\mathbb{E}_{t_k} ( \{ Y_{t_{k+1}} + \int_{t_k}^{t_{k+1}} f(\theta_s) ds \} \Delta
W_k^*)$, with a bounded $f$, it follows that 
\begin{equation}\label{eqs7}
  h^2 |\overline{Z}_{t_k}-Z_{t_k}^N|^2 \leq d\ h\left(
    \mathbb{E}_{t_k}(|\Delta Y_{k+1}|^2)-
    |\mathbb{E}_{t_k}(\Delta Y_{k+1})|^2\right)+Ch^3.                                                
\end{equation}
By plugging (\ref{eqs6}) and (\ref{eqs7}) into (\ref{eqs5}), we get 
\begin{align*}
  |\Delta Y_k|^2 & \leq  \displaystyle (1+ \gamma h)(\mathbb{E}_{t_k}(\Delta Y_{k+1}))^2\\
  &  + \displaystyle C( h+
  \frac{1}{\gamma }) [h^2 +  \mathbb{E}_{t_k}\int_{t_k}^{t_{k+1}}
  |X_s-X^N_{t_k}|^2 ds + \mathbb{E}_{t_k}\int_{t_k}^{t_{k+1}} |Y_s-Y^N_{t_{k+1}}|^2 ds]\\
  &  + \displaystyle  C( h+
  \frac{1}{\gamma })[\mathbb{E}_{t_k}   \int_{t_k}^{t_{k+1}} |Z_s-\overline{Z}_{t_k}|^2 ds
  +  \mathbb{E}_{t_k}(|\Delta Y_{k+1}|^2)-
  |\mathbb{E}_{t_k}(\Delta Y_{k+1})|^2  ].
\end{align*}
We can write $\E_{t_k} |Y_s-Y^N_{t_{k+1}}|^2 \leq 2\E_{t_k}|Y_s-Y_{t_{k+1}}|^2
+2\E_{t_k}|\Delta Y_{k+1}|^2$. By doing the same for $X_s-X^N_{t_{k+1}}$, and taking
$\gamma=C$, we obtain
\begin{align*}
  |\Delta Y_k|^2  \leq &\displaystyle (1+ Ch)\mathbb{E}_{t_k}|\Delta
  Y_{k+1}|^2+Ch|\Delta X_k|^2+Ch\mathbb{E}_{t_k}|\Delta Y_{k+1}|^2\\
  & +\displaystyle C [h^2 +  \mathbb{E}_{t_k}\int_{t_k}^{t_{k+1}}
  |X_s-X_{t_k}|^2 ds + \mathbb{E}_{t_k}\int_{t_k}^{t_{k+1}} |Y_s-Y_{t_{k+1}}|^2 ds]\\
  &+\displaystyle  C [\mathbb{E}_{t_k}   \int_{t_k}^{t_{k+1}}
  |Z_s-\overline{Z}_{t_k}|^2 ds].
\end{align*}
Using (\ref{eqs18}) yields $|\Delta Y_k|^2 \leq \displaystyle (1+ Ch)\mathbb{E}_{t_k}|\Delta
Y_{k+1}|^2+Ch|\Delta X_k|^2 +Ch^2.$\qed
\vspace{-0.5cm}
\subsection{Proof of  $\mathbb{E}
  (\sum_{k=0}^{N-1} \int_{t_k}^{t_{k+1}} |Z^N_{t_k}-Z_t|^2 dt)^{\frac{p}{2}}=
  O(h^p).$}
\vspace{-0.5cm}
First of all, we can split this summation into two terms 
\begin{equation*}
  \mathbb{E}\big( \sum_{k=0}^{N-1} \int_{t_k}^{t_{k+1}} |Z^N_{t_k} - Z_t|^2
  dt\big)^p\leq  C\mathbb{E}\big( \sum_{k=0}^{N-1} \int_{t_k}^{t_{k+1}}
  |\overline{Z}_{t_k} - Z_t|^2 dt\big)^p + C\mathbb{E}\big( h\sum_{k=0}^{N-1}
  |\Delta \overline{Z}_k|^2\big)^p.
\end{equation*}
Thanks to (\ref{eqs18}), we have $\mathbb{E}\big( \sum_{k=0}^{N-1} \int_{t_k}^{t_{k+1}}
  |\overline{Z}_{t_k} - Z_t|^2 dt\big)^p \leq T^{p-1}\sum_{k=0}^{N-1}\int_{t_k}^{t_{k+1}}
\mathbb{E}|\overline{Z}_{t_k} - Z_t|^{2p} dt=O(h^p).$

{\it Scheme of the proof of $\mathbb{E}\big( h\sum_{k=0}^{N-1}
  |\Delta \overline{Z}_k|^2\big)^p=O(h^p)$. } The first key point is to slice the summation into small intervals and show that the result
is true for small time intervals. The second key point is to use Rosenthal's
inequality, see Theorem $2.12$ page $23$ of \cite{hall_80}. By using (\ref{eqs7}) and taking the
expectation, we can write :
\begin{equation}\label{eqs8}
  \mathbb{E}\big(h \sum_{k=0}^{k_1}|\Delta \overline{Z}_k|^2\big)^p \leq \displaystyle C \mathbb{E}  \big(
  \sum_{k=0}^{k_1} \mbox{Var}_{t_k} \Delta Y_{k+1} \big)^p+ C h^p.
\end{equation}
We use Rosenthal's inequality to upper bound
\begin{align*}
  \mathbb{E}  \big(
  \sum_{k=0}^{k_1} \mbox{Var}_{t_k} \Delta Y_{k+1} \big)^p &\leq \displaystyle  C \mathbb{E}\big( \sum_{k=0}^{k_1} \Delta Y_{k+1} -\mathbb{E}_{t_k} \Delta Y_{k+1}
  \big)^{2p},\\
  &\leq  \displaystyle C 3^{2p-1}\big[ \mathbb{E} \Delta Y_{k_1+1}^{2p} + \mathbb{E} \Delta
  Y_0^{2p} + \mathbb{E}\big(\sum_{k=0}^{k_1}(\Delta Y_k
  -\mathbb{E}_{t_k} \Delta Y_{k+1})\big)^{2p}\big].
\end{align*}
By plugging this inequality into (\ref{eqs8}) and using the previous estimate on $|\Delta Y_k|$, we
get
\begin{equation}\label{eqs9}  
  \mathbb{E}\big(h \sum_{k=0}^{k_1}|\Delta \overline{Z}_k|^2\big)^p \leq  O(h^p) +C \mathbb{E}\big(\sum_{k=0}^{k_1}(\Delta Y_k
  -\mathbb{E}_{t_k} \Delta Y_{k+1})\big)^{2p}.
\end{equation}
We now tackle the term $\Delta Y_k-\mathbb{E}_{t_k} \Delta Y_{k+1}$. Using
(\ref{eqs0}), we have $\sum_{k=0}^{k_1} (\Delta Y_k -  \mathbb{E}_{t_k} \Delta Y_{k+1}=\sum_{k=0}^{k_1}\int_{t_k}^{t_{k+1}} (\mathbb{E}_{t_k} (f^N_{t_k}-f(\theta_s))
  ) \mbox{ds}.$
By doing the same kind of proof as before, that is using the fact that $f$ is Lipschitz
and the results on $\E|\Delta X_k|^{2p}$ and $\E|\Delta Y_k|^{2p}$, we find
\begin{equation*}
  \displaystyle \mathbb{E}\big(\sum_{k=0}^{k_1} (\Delta Y_k -  \mathbb{E}_{t_k}
  \Delta Y_{k+1})\big)^{2p}  \leq  O(h^p) + C (hk_1)^p\mathbb{E} \big(h \sum_{k=0}^{k_1}
  |\Delta \overline{Z}_k|^2\big)^p. 
\end{equation*}
By plugging this term back into (\ref{eqs9}), we can write $\big(1-C (hk_1)^p\big)\mathbb{E}\big(h \sum_{k=0}^{k_1}|\Delta \overline{Z}_k|^2\big)^p =O(h^p).$
Consequently, if we choose $k_1 \leq  \frac{1}{(2C)^{\frac{1}{p}}h}$ we come up with
  $\mathbb{E}\big(h \sum_{k=0}^{k_1}|\Delta \overline{Z}_k|^2\big)^p = O(h^p).$
This result can be extended to any summation involving at most $\Delta k$ terms, where  $\Delta k \leq
\frac{1}{(2 C)^{\frac{1}{p}} h}$. We can cover the interval $\{0,\cdots,N-1\}$ with a finite
number of elementary intervals of size $\Delta k$ and we get $\mathbb{E}\big(h
\sum_{k=0}^{N-1}|\Delta \overline{Z}_k|^2\big)^p = O(h^p)$, which completes our
proof. \qed\\
From this result and (\ref{eqs18}), we also deduce
\begin{equation}\label{eqs10}
  \mathbb{E}\big(h \sum_{k=0}^{N-1}|\Delta Z_k|^2\big)^p = O(h^p),
\end{equation}
which is very useful in the following.

\section{Proof of Theorem \ref{theo2}.}
\label{section:proof:theo2}
To expand the error, we use usual techniques of stochastic analysis, combining martingale estimates and Malliavin calculus tools.

\subsection{Preliminary estimates}
Sections \ref{section:proof:theo2} and \ref{section:proof:theo4} contain proofs
with similar calculations, which are quite technical. In order to be as clear as possible, we state two results
really useful in the sequel, which are related to Malliavin calculus (see \cite{nualart_95}). The results give sufficient conditions for expectations and conditional expectations to be small w.r.t. the time step $h$. They are based on ideas from \cite{koha:pett:02} and \cite{munos_04}.

\begin{prop}\label{prop5}
  Let $F \in \D^{1,2}$ with $\E_{t_k} |F|^2 + \sup_{t_k \le s \leq T} \E_{t_k} |D_s F|^2 < \infty$ and let $U$ be an Itô process of the form $U_t=U_0+ \int_0^t \alpha_s ds +
  \int_0^t \beta_s dW_s$, with $\sup_{t_k \leq s\le T}\E_{t_k}|\alpha_s|^2 + \sup_{t_k \leq s\le T}\E_{t_k}|\beta_s|^2< \infty$. Then, $\forall (t,t')$ such that $t_k \leq t \leq
  t' \leq t_{k+1}$,
  \begin{align*}
    |\E_{t_k}[F(U_t -U_{t'})]|\leq (t'-t)\big[&(\E_{t_k}
  |F|^2)^{\frac{1}{2}}(\sup_{t \leq s \leq t'} \E_{t_k} |\alpha_s|^2)^{\frac{1}{2}}\\
    &+(\sup_{t \leq s \leq t'}\E_{t_k}
    |D_sF|^2)^{\frac{1}{2}}(\sup_{t \leq s \leq t'} \E_{t_k}
    |\beta_s|^2)^{\frac{1}{2}}\big].
    \end{align*}
\end{prop}
This proposition can be easily proved. Assume without loss of generality that $F$ and $U$ are one-dimensional. From the duality formula, we have $\E_{t_k}[F(\int_t^{t'} \alpha_s ds+\int_t^{t'} \beta_s
dW_s)]=\E_{t_k}[\int_t^{t'} (F\alpha_s+D_s F\cdot \beta_s) ds]$. Thanks to Cauchy Schwarz inequality and
hypotheses on $\alpha$ and $\beta$, we get the result.
\begin{defn} $F$ satisfies the condition $R_k$ if $F \in
  \mathbb{D}^{k,\infty}$ and if $\mathcal{C}_{k,p}(F):=\|F\|_{L_p}+\sum_{j \le k}\sup_{0 \leq
    s_1,...,s_j\leq T} \|D_{s_1,..,s_j} F\|_{L_p} < \infty$.
  \end{defn}
\begin{prop}\label{prop6} Let $F$ satisfy the condition $R_2$.
For simplicity we set $dW^0_s=ds$.
Assume that $U_t \in \mathbb{R}^d$
satisfies the following stochastic expansion property
\begin{equation}\label{eq:se}\tag{\text{$\mathcal{P}$}}
  U_t=\sum_{i,j=0}^{q} c_{i,j}^{U,0}(t)\int_0^t c_{i,j}^{U,1}(s)
  \big(\int_{\eta(s)}^s c_{i,j}^{U,2}(r)dW^i_r\big)dW^j_s,
\end{equation}
where $\{(c_{i,j}^{U,i_1}(t))_{t \geq 0}$ : $0\leq i,j \leq q, 0\leq i_1 \leq 2\}$
are adapted processes satisfying
\begin{itemize}
\item $\forall (i,j), 1\leq i,j \leq q$, $\forall t \in [0,T]$, $c^{U,0}_{i,j}(t)$
  satisfies $R_2$, and \linebreak[4]$\mathcal{C}^U_{2,p}:=\sup_{0 \le t \le T}\sup_{1 \le i,j \le q}
  \mathcal{C}_{2,p}(c_{i,j}^{U,0}(t))< \infty$, $p\ge 1$.\\
\item $\forall (i,j), 1\leq i,j \leq q$, $\forall t \in [0,T]$, $c^{U,1}_{i,j}(t),c^{U,0}_{0,j}(t),c^{U,0}_{i,0}(t),c^{U,1}_{i,0}(t)$
  satisfy $R_1$, and \linebreak[4] $\mathcal{C}^U_{1,p}:=\sup_{0 \le t \le T}\sup_{1 \le i,j \le q}
  \{\mathcal{C}_{1,p}(c_{i,j}^{U,1}(t))+\mathcal{C}_{1,p}(c_{0,j}^{U,0}(t))+\mathcal{C}_{1,p}(c_{i,0}^{U,0}(t))+\mathcal{C}_{1,p}(c_{i,0}^{U,1}(t))\}\break < \infty$, $p\ge 1$.\\
\item $\forall (i,j), 0\leq i,j \leq q$, $\forall t \in [0,T]$,
  $c^{U,2}_{i,j}(t),c^{U,1}_{0,j}(t), c^{U,0}_{0,0}(t)$
  satisfy $R_0$, and $\mathcal{C}^U_{0,p}:=\sup_{0 \le t \le T}\sup_{0 \le i,j \le q}
   \{\mathcal{C}_{0,p}(c_{i,j}^{U,2}(t))+\mathcal{C}_{0,p}(c_{0,j}^{U,1}(t))+\mathcal{C}_{0,p}(c_{0,0}^{U,0}(t))\} < \infty$, $p\ge 1$.
\end{itemize}
Thus, there is a constant $K(T)$ which depends polynomially on $\mathcal{C}_{2,p}(F), \mathcal{C}^U_{2,p},
\mathcal{C}^U_{1,p},\mathcal{C}^U_{0,p}$ (for some $p \ge 1$) such that $|\E[FU_t]|\leq K(T)h  .$
\end{prop}
Indeed, we have
\begin{align*}
  \E(FU_t)
  =&\sum_{i,j=0}^q \E\big(F c_{i,j}^{U,0}(t)\int_0^t c_{i,j}^{U,1}(s)(\int_{\eta(s)}^s   c_{i,j}^{U,2}(r)dW_r^i)dW_s^j\big)\\
  =&\sum_{i,j=1}^q \int_0^t \int_{\eta(s)}^s \E \big( D^i_r \big[ D^j_s \{F
  c_{i,j}^{U,0}(t)\} c_{i,j}^{U,1}(s)\big]c_{i,j}^{U,2}(r)\big)\mbox{dr ds}\\
  &+\sum_{j=1}^q \int_0^t \int_{\eta(s)}^s \E  \big[ D^j_s \{F
  c_{0,j}^{U,0}(t)\} c_{0,j}^{U,1}(s)c_{0,j}^{U,2}(r)\big]\mbox{dr ds}\\
  &+\sum_{i=1}^q \int_0^t \int_{\eta(s)}^s \E  \big[ D^i_r \{F
  c_{i,0}^{U,0}(t) c_{i,0}^{U,1}(s)\}c_{i,0}^{U,2}(r)\big]\mbox{dr ds}\\
  &+\int_0^t \int_{\eta(s)}^s \E\big( F
  c_{0,0}^{U,0}(t) c_{0,0}^{U,1}(s)c_{0,0}^{U,2}(r)\big)\mbox{dr ds}.
\end{align*}
Then, the result readily follows.
  \begin{rem}\label{remarque4}Under Hypothesis \ref{hyp4}, we can show (see later the proof of \eqref{eqs37}) that for each $t$, $X^N_t -X_t$ satisfies the expansion \ref{eq:se}. Hence, if $F$ satisfies $R_2$, Proposition \ref{prop6} yields $$|\E[F(X^N_t -X_t)]|=O(h)$$ uniformly in $t\in[0,T]$,
which is a very useful result for the sequel.
\end{rem}
\subsection{Expansion of $Y^N_{t_k}-Y_{t_k}$}
In the following, we assume that Hypothesis \ref{hyp4} is in force. This implies in
particular that $u$ is bounded, of class $C_b^{3/2,3}$ (see
Theorem \ref{theo5}). We also easily prove that $\forall p \geq 1, \forall k \in
\{0,\cdots,N-1\}$ (see \cite{nualart_95} e.g.)
\begin{align}\notag
 \bullet\ &\E_{t_k}(\sup_{t_k \le t \le T}|X_t|^{2p})<K(T)(1+|X_{t_k}|^{2p}), \sup_{t_k \le s
    \leq T} \E_{t_k}(\sup_{t_k \le t \le T}|D_s X_t|^p) \le C,\\
  &\sup_{t_k \le s,r
   \leq T}
 \E_{t_k}(\sup_{t_k \le t \le T}|D_r D_sX_t|^p)+\sup_{t_k \le s,r,v
   \leq T}
 \E_{t_k}(\sup_{t_k \le t \le T}|D_v D_r D_sX_t|^p) \leq C,\label{eqs48}  \\
\bullet\ &\E_{t_k}(\sup_{t_k \le t \le T}|X^N_t|^{2p})<K(T)(1+|X^N_{t_k}|^{2p}),\sup_{N,t_k \le
   s \leq T} \E_{t_k}(\sup_{t_k \le t \le T}|D_s X^N_t|^p)\le C,\notag\\ &
 \sup_{N,t_k \le s,r
   \leq T}
 \E_{t_k}(\sup_{t_k \le t \le T} |D_r D_sX^N_t|^p)+\sup_{N,t_k \le s,r,v
   \leq T}
 \E_{t_k}(\sup_{t_k \le t \le T} |D_v D_r D_sX^N_t|^p) \leq C \label{eqs49}.
\end{align}
Due to the Markov property of $(X^N_{t_k})_k$, one has $Y^N_{t_k}=u^N(t_k,X^N_{t_k})$ for some Lipschitz function $u^N(t_k,\cdot)$ (see \cite{lemor_05}) with an obvious definition of $u^N$. Actually, under our assumptions, this function is even three times differentiable w.r.t. $x$. Thus, the difference $\Delta Y_k$ can be written as follows:
\begin{align*}
  \Delta Y_k=(u^N(t_k,X^N_{t_k})-u(t_k,X^N_{t_k}))+(u(t_k,X^N_{t_k})-u(t_k,X_{t_k})).
\end{align*}
Since $u$ is of class $C_b^{3/2,3}$, the last term of the previous inequality becomes
\begin{equation}\label{eqs13}
  u(t_k,X^N_{t_k})-u(t_k,X_{t_k})=\nabla_x u(t_k,X_{t_k}) \Delta X_k +O(|\Delta
  X_k|^2).
\end{equation}
To complete the proof, we apply the following lemma
\begin{lem}\label{lemme1} Under Hypothesis \ref{hyp4}, $ |u^N(t_k,x)-u(t_k,x)|\leq K(T,x)h$.
  \end{lem}
The result above is new but not so surprising. Indeed, if $f$ is identically zero, the difference is only related to the weak approximation of $\Phi(X_T)$ by $\Phi(X^N_T)$: from \cite{bally_96}, one knows that this is of order $h$.

The rest of this section is devoted to the proof of the lemma. We only give the proof for $t_k=0$. We want to find a upper bound for $|u^N(0,x)-u(0,x)|=|\Delta Y_0|$.

For the sake of clarity, we split the proof into several steps.

{\bf Step 1 : linearization of the error.} We show that 
\begin{equation}
  \label{eq:y:dl}
  \Delta Y_k= \E_{t_k}(\Delta Y_{k+1} \xi_k + h f'_x(\theta_{t_k}) \Delta X_k + h \chi_k),
\end{equation}
  with 
  \begin{align}
 &   \xi_k =(1+hf_y'(\theta_{t_k}) + f_z'(\theta_{t_k})\Delta W_k),\label{eq:xi}\\
\label{eq:chi}&\chi_k=\int_{t_k}^{t_{k+1}} (G_0(s,X_s)+  f_y'(\theta_{t_k})G_y(s,X_s)+ f_z'(\theta_{t_k})G_z(s,X_s))ds\\
\nonumber&+\int_0^1 (1-\lambda)\big[ \Delta X_k^* f^{''}_{xx}(\theta^{\lambda}_{t_k})\Delta
    X_k+f^{''}_{yy}(\theta^{\lambda}_{t_k})(Y^N_{t_{k+1}}-Y_{t_k})^2+\Delta Z_k
    f^{''}_{zz}(\theta^{\lambda}_{t_k})\Delta Z_k^*\\
&     +2\Delta X_k^*f^{''}_{xy}(\theta^{\lambda}_{t_k})(Y^N_{t_{k+1}}-Y_{t_k})+2\Delta X_k^*f^{''}_{xz}(\theta^{\lambda}_{t_k})\Delta Z_k^*+ 2(Y^N_{t_{k+1}}-Y_{t_k})f^{''}_{yz}(\theta^{\lambda}_{t_k})\Delta
    Z_k^*\big]d\lambda,
\nonumber  \end{align}
    where $\theta^{\lambda}_{t_k}=\lambda(t_k,X^N_{t_k},Y^N_{t_{k+1}},Z^N_{t_k})+(1-\lambda)\theta_{t_k}$ and $G_0,G_y,G_z$ are bounded functions.
From (\ref{eqs0}) and by introducing $f(\theta_{t_k})$, we have
\begin{align}\label{eqs52}
  \Delta Y_k= \E_{t_k}\big(\Delta Y_{k+1}+h(f^N_{t_k}-f(\theta_{t_k})) +\int_{t_k}^{t_{k+1}} ( f(\theta_{t_k})-f(\theta_s)) \mbox{ds}\big).
\end{align} 
By applying Itô's formula to $f(\theta_u)$ between $t_k$ and $s$ we show that, under
Hypothesis \ref{hyp4} , $\int_{t_k}^{t_{k+1}} \E_{t_k}(f(\theta_{t_k})-f(\theta_s))ds= h \int_{t_k}^{t_{k+1}}E_{t_k}(G_0(s,X_s))ds$, where $G_0$ is a bounded function. In the second
term, perform a second order expansion of $f$ around $\theta_{t_k}$ to get
 \begin{align}\label{eqs11}
    &f^N_{t_k}-f(\theta_{t_k}) =f_x'(\theta_{t_k})\Delta X_k+f_y'(\theta_{t_k})\Delta  Y_{k+1}+f_z'(\theta_{t_k})\Delta
    Z_k^*+ f'_y(\theta_{t_k})(Y_{t_{k+1}}-Y_{t_k})\notag\\
    &+\int_0^1 (1-\lambda)\big[ \Delta X_k^* f^{''}_{xx}(\theta^{\lambda}_{t_k})\Delta
    X_k+f^{''}_{yy}(\theta^{\lambda}_{t_k})(Y^N_{t_{k+1}}-Y_{t_k})^2+\Delta Z_k
    f^{''}_{zz}(\theta^{\lambda}_{t_k})\Delta Z_k^*\\
&+2\Delta
    X_k^*f^{''}_{xy}(\theta^{\lambda}_{t_k})(Y^N_{t_{k+1}}-Y_{t_k})+2\Delta X_k^*f^{''}_{xz}(\theta^{\lambda}_{t_k})\Delta Z_k^*+ 2(Y^N_{t_{k+1}}-Y_{t_k})f^{''}_{yz}(\theta^{\lambda}_{t_k})\Delta
    Z_k^*\big]d\lambda. \notag
 \end{align}
Note that $\Ek(Y_{t_{k+1}}-Y_{t_k})=\Ek\int_{t_k}^{t_{k+1}} G_y(s,X_s)ds$. If we closely look at (\ref{eqs11}), we can see that we need to develop $\Delta
Z_k$. By using (\ref{eqs3}), we can write 
\begin{align*}
Z^N_{t_k}=\frac{1}{h}\E_{t_k}(\Delta Y_{k+1} \Delta W_k^*)+\frac{1}{h}\E_{t_k}(u(t_{k+1},X_{t_{k+1}})\Delta W_k^*).
\end{align*}
Introducing the weak derivative of $X_{t_{k+1}}$ (see \cite{nualart_95} p.109), the second term of this summation equals
$\frac{1}{h}\E_{t_k} \int_{t_k}^{t_{k+1}} \nabla_x u(t_{k+1},X_{t_{k+1}}) D_t
X_{t_{k+1}}\mbox{dt}$, where $D_t X_{t_{k+1}}= \nabla_x X_{t_{k+1}} (\nabla_x X_t)^{-1}
  \sigma(t,X_t)$. Since $Z_{t_k}=\nabla_x  u(t_k,X_{t_k}) \sigma(t_k,X_{t_k})$, one gets
\begin{align*}
  &\Delta Z_k=\frac{1}{h}\E_{t_k}(\Delta Y_{k+1}\Delta W_k^*)\\
  &+\frac{1}{h}\int_{t_k}^{t_{k+1}}\E_{t_k}\big(\nabla_x
  u(t_{k+1},X_{t_{k+1}})\nabla_x X_{t_{k+1}} (\nabla_x X_t)^{-1}
  \sigma(t,X_t)
  -\nabla_x  u(t_k,X_{t_k}) \sigma(t_k,X_{t_k})\big) \mbox{dt}.
  \end{align*}
The term in the second conditional expectation is equal to $\nabla_x
  u(t_{k+1},X_{t_{k+1}}) \break \nabla_x X_{t_{k+1}} (\nabla_x X_t)^{-1}
  \sigma(t,X_t)\pm \nabla_x  u(t,X_{t}) \sigma(t,X_{t})
  -\nabla_x  u(t_k,X_{t_k}) \sigma(t_k,X_{t_k})$: hence, two applications of
  Itô's formula (for the first contribution between $t$ and $t_{k+1}$, for the second one between $t_k$ and $t$) prove that
\begin{equation}\label{eqs14}
  \Delta Z_k^*=\int_{t_k}^{t_{k+1}}\E_{t_k}(G_z(s,X_s)) \mbox{ds}+\frac{1}{h}\E_{t_k}(\Delta Y_{k+1}\Delta W_k),
  \end{equation}
for a bounded function $G_z$. Plugging this equality and (\ref{eqs11}) into
(\ref{eqs52}) yields \eqref{eq:y:dl}.

{\bf Step 2 : another formula of $\Delta Y_0$.}
First of all, we replace $Y^N_{t_{k+1}}-Y_{t_k}$ by $\Delta Y_{k+1}+
 Y_{t_{k+1}}-Y_{t_k}$ in the expression of $\chi_k$. Then, easy computations combining Proposition \ref{prop2} and estimates (\ref{eqs18}) show that 
 \begin{equation}
   \label{eq:chik}
   \tilde{\chi}_k=\Ek(\chi_k)=O_k(h) + O(|\Delta X_k|^2+ |\Delta Z_k|^2).
 \end{equation}
From \eqref{eq:y:dl}, we
 deduce the following equality
\begin{align}\label{eqs12}
  \Delta Y_0=\E(\Delta Y_N\xi_0\cdots\xi_{N-1} + h \sum_{i=0}^{N-1} (f'_x(\theta_{t_i}) \Delta X_i+\tilde{\chi}_i)
  \xi_0 \cdots\xi_{i-1}).
\end{align}
Now it is enough to show that all terms of this summation are $O(h)$.
In the following, $\eta_0=1$ and $\eta_i=\xi_0\cdots\xi_{i-1}$ for $i\le N$.

{\bf Step 3 : some results on $\eta_N=\xi_0\cdots\xi_{N-1}$.}\\
We establish the following results on $\eta_N$:
\begin{align}
  & \eta_k \text{ satisfies the condition } R_2 \text{ uniformly in } k,  \mbox{ i.e.
  } \forall k, \eta_k \in \mathbb{D}^{2,\infty}\notag\\
  \label{eq:a}& \text{ and }\max_{k \le N}\mathcal{C}_{2,p}(\eta_k)<\infty, \forall p \ge 1,\\
 \label{eq:b}&\E(\max_{0 \leq k \leq N} |\eta_k|^p)+\sup_{r \le T}\E(\max_{0 \le k \le
   N} |D_r \eta_k|^p)+\sup_{r,s \le T}\E(\max_{0 \le k \le N} |D_rD_s \eta_k|^p)< \infty.
\end{align}
{\it Proof of \eqref{eq:a}.} We have $\eta_0=1$, and for $i\ge 1$
\begin{equation}\label{eqs15}
   \eta_i=\eta_{i-1}(1+hf_y'(\theta_{t_{i-1}})+f_z'(\theta_{t_{i-1}})\Delta W_{i-1}).
 \end{equation}
 We begin to show that $\max_{k\le N} \|\eta_k\|_{L_p}=O(1)$ for $p\geq 1$. Since $f'_y$ and $f'_z$ are bounded, we easily prove that $\E_{t_{i-1}}(1+hf_y'(\theta_{t_{i-1}})+f_z'(\theta_{t_{i-1}})\Delta W_{i-1})^{2p}\leq (1+C h)$, whence $\E|\eta_i|^{2p}\le (1+Ch)\E|\eta_{i-1}|^{2p}.$
    We deduce that $\max_{k\le N} \|\eta_k\|_{L_p} =O(1)$. 

Now, let us show that $\max_{k \le N}\E|D_r \eta_k|^p=O(1)$, uniformly in $r$. Let $r$ be such that
$t_{k-1} < r \leq t_k$. $\forall i \leq k-1, D_r \eta_i=0$. We note that $D_r \eta_k=\eta_{k-1}f_z'(\theta_{t_{k-1}})$. For $i \geq k+1$, we have
   \begin{align}\nonumber
     D_r \eta_i&= D_r
     \eta_{i-1}+hD_r(\eta_{i-1}f_y'(\theta_{t_{i-1}}))+\sum_{l=1}^q D_r(\eta_{i-1}f'_{z_l}(\theta_{t_{i-1}}))\Delta W^l_{i-1},\\
     &=\eta_{k-1}f_z'(\theta_{t_{k-1}}) +h\sum_{j=k}^{i-1}D_r(\eta_j
     f_y'(\theta_{t_j}))+\sum_{l=1}^q \sum_{j=k}^{i-1}D_r(\eta_jf'_{z_l}(\theta_{t_j}))\Delta W^l_j.\label{eq:c1}
     \end{align}
Applying Burkholder-Davis-Gundy's inequality to the martingale
\linebreak[4] $\sum_{j=k}^{i-1} D_r(\eta_j f'_{z_l}(\theta_{t_j})) \Delta W^l_j$ yields
\begin{align*}
  \E|D_r \eta_i|^p& \leq C\E|\eta_{k-1}|^p+C_ph\sum_{j=k}^{i-1}\E|D_r(\eta_j
     f_y'(\theta_{t_j}))|^p+C\sum_{l=1}^q \E|h\sum_{j=k}^{i-1}|D_r(\eta_jf'_{z_l}(\theta_{t_j}))|^2|^{\frac{p}{2}}\\
     & \leq C\E|\eta_{k-1}|^p+Ch\sum_{j=k}^{i-1}\E|D_r(\eta_j
     f_y'(\theta_{t_j}))|^p+C \sum_{l=1}^q h\sum_{j=k}^{i-1}\E|D_r(\eta_j
     f'_{z_l}(\theta_{t_j}))|^p\\
     & \leq C(1+\E|\eta_{k-1}|^p)+Ch\sum_{j=k+1}^{i-1}\E|D_r\eta_j|^p,
     \end{align*}
using the boundedness of the derivatives of $f$, $\max_{j\le N}\|\eta_j\|_q =O(1)$, idendity (\ref{eqs41}), $u,\sigma \in
C^{1,2}_b$, and estimates (\ref{eqs48}).
By applying Gronwall's lemma, we get $\max_{k\leq i \leq N}\E|D_r\eta_i|^p \leq C
(1+\E|\eta_{k-1}|^p), t_{k-1} < r \leq t_k$.\\
 Then, $\max_{k \le N}\E|D_r \eta_k|^p=O(1)$, uniformly in $
r \in [0,T]$.
The proof concerning the derivative of order $2$ can be done following the same
scheme.\qed

{\it Proof of \eqref{eq:b}}. We begin to show that $\E(\max_{k\leq N}|\eta_k|^p)< \infty$. The idea is to use a martingale property in order to apply Doob's inequality.
Since $\eta_i=\eta_{i-1} + h \eta_{i-1} f_y'(\theta_{t_{i-1}})+ \eta_{i-1} f_z'(\theta_{t_{i-1}})\Delta W_{i-1}$, one has $\eta_k=1+\sum_{i=1}^{k} h \eta_{i-1} f_y'(\theta_{t_{i-1}})+ \eta_{i-1} f_z'(\theta_{t_{i-1}})\Delta W_{i-1}$. Thus, 
\begin{align*}
   \E(\max_{k\leq N}|\eta_k|^p) \leq C \big(1+ \E(\sum_{i=1}^{N} h |\eta_{i-1}| |f_y'(\theta_{t_{i-1}})|)^p +\E(\max_{k\leq N}|\sum_{i=1}^{k} \eta_{i-1} f_z'(\theta_{t_{i-1}})\Delta W_{i-1}|^p)\big).
\end{align*}
The last term is upper bounded by $C \E(h \sum_{i=1}^{N} |\eta_{i-1}
f_z'(\theta_{t_{i-1}})|^2)^{\frac{p}{2}}\leq C  h\break \sum_{i=1}^{N} \E |\eta_{i-1}
f_z'(\theta_{t_{i-1}})|^p$. Using the estimate \eqref{eq:a}, we get  $\E(\max_{k\leq
  N}|\eta_k|^p)< \infty$.\\
To prove that $\sup_{r \le T}\E(\max_{k\leq N}|D_r\eta_k|^p)< \infty$, we proceed in
the same way, by starting from (\ref{eq:c1}). For the second derivative, this is analogous.

{\bf Step 4 : we prove that $\E(\Delta Y_N\eta_N)=O(h)$.}\\
If $\eta_N$ were equal to 1, the results of \cite{bally_96} would directly apply. Here the approach has to be different and we use techniques of Malliavin calculus. We have $\E(\Delta
Y_N\eta_N)=\E(\eta_N\Phi(X^N_T)-\eta_N\Phi(X_T))$. Let us introduce $X^{N,\lambda }_t=(1-\lambda ) X_t
+\lambda X^N_t$. Thus, we have
\begin{equation*}
\E(\Delta Y_N\eta_N)=\int_0^1 \E\big(\eta_N
\Phi'_x(X^{N,\lambda }_T)(X^N_T-X_T)\big)\mbox{d}\lambda. 
\end{equation*}
As $\Phi \in
C^{3+\alpha}$, by using (\ref{eq:a}), (\ref{eqs48}) and (\ref{eqs49}),
 we note that $\eta_N\Phi'_x(X^{N,\lambda }_T)$ satisfies $R_2$. By applying Remark \ref{remarque4}, we deduce that $\E(\Delta Y_N\eta_N)=O(h).$
 
{\bf Step 5 : we prove that $\E(f'_x(\theta_{t_i})\Delta X_i \eta_i)=O(h)$.} This is a very similar proof to Step 4, in a case where $\Phi(x)=x$.

{\bf Conclusion.} We now work on $h\E(\sum_{i=0}^{N-1} \tilde{\chi}_i \eta_i)$, where $|\tilde{\chi}_k|\le \lambda^N_k h + K(T,x)|\Delta X_k|^2+K(T,x) |\Delta Z_k|^2$. Hence, 
\begin{align*}
 | h \sum_{i=0}^{N-1}\E(\tilde{\chi}_i \eta_i)|\leq &C\sum_{i=0}^{N-1}\E(
  \lambda^N_i |\eta_i|)h^2+K(T,x)\sum_{i=0}^{N-1}h\E\big( |\eta_i|(|\Delta X_i|^2
  +|\Delta Z_i|^2)\big)\\
  \leq &K(T,x)h+K(T,x)\sum_{i=0}^{N-1}h\E( |\eta_i| |\Delta Z_i|^2)\\
\leq & K(T,x)h+ K(T,x)\big( \E(\max_{0 \leq i \leq N-1}
    |\eta_i|)^2\big)^{\frac{1}{2}}\big(\E(h\sum_{i=0}^{N-1} |\Delta
  Z_i|^2)^2\big)^{\frac{1}{2}}.
\end{align*}
By using \eqref{eq:b} on $(\eta_i)_i$ and the upper bound \eqref{eqs10} we get that
 $|h\E(\sum_{i=0}^{N-1} \tilde{\chi}_i \eta_i)|\leq K(T,x)h.$
By combining this result and the results of Step 4 and Step 5, (\ref{eqs12}) shows that
$|\Delta Y_0|\leq K(T,x)h.$ Lemma \ref{lemme1} is proved. \qed

\section{Proof of Theorem \ref{theo4}.}
\label{section:proof:theo4}
As it could be expected, its proof is more difficult. The main extra ingredient is the convergence of the weak derivative of the discrete BSDE $(Y^N,Z^N)$, with the rate of convergence $N^{-1/2}$. The next paragraph is aimed at proving this result. In the following, Hypothesis \ref{hyp5} is in force. 
\subsection{ Proof of an intermediate result}
\begin{prop}\label{lemme3}
  Let $r \in ]0,t_1[$. Under Hypothesis \ref{hyp5}, we have
$    \max_{1 \leq i \leq N}\E|D_r \Delta Y_{i}|^2+h\E\big(\sum_{i=1}^{N-1}|D_r
    \Delta Z_i^*|^2 \big)=O(h)$, uniformly in $r$.
\end{prop}
This proposition is analogous to Theorem \ref{theo1}, where $q=2$, and the scheme of its proof as well. However, there is a significative difference: the BSDE solved by the weak derivatives (see (\ref{eqs22}-\ref{eqs23}-\ref{eqs24})) has a non Lipschitz driver, which requires extra technicalities that we detail. In what follows, we fix $r \in ]0,t_1[$ and introduce some specific notations. $\widehat{X_t}$ stands for $D_r X_t$. In the case of $Z_t$, which is a row vector, $\widehat{Z_t}$ is a matrix whose the $i$-th column is $D^i_r Z_t^*$. It is well-known (Proposition $5.3$ of \cite{karoui_97}) that $(\widehat{Y_t}, \widehat{Z_t})_{r\le t\le T}$ solves
\begin{align}\label{eqs22}
  \widehat{Y_t}= \Phi'_x(X_T)\widehat{X_T}+\int_t^T  (f'_x(\theta_s)
  \widehat{X_s} +  f'_y(\theta_s)\widehat{Y_s} + f'_z(\theta_s)\widehat{Z_s}) ds- (\int_t^T {\widehat{Z_s}}^* dW_s)^*.
\end{align}
Regarding $(\widehat{Y^N},\widehat{Z^N})$, one obtains
\begin{align}\label{eqs23}
  \widehat{Y^N_{t_k}}=&\E_{t_k}[\widehat{Y^N_{t_{k+1}}}+h \nabla_x
  f^N_{t_k}\widehat{X^N_{t_k}}+h\nabla_y
  f^N_{t_k}\widehat{Y^N_{t_{k+1}}}+h\nabla_z
  f^N_{t_k}\widehat{Z^N_{t_k}}],\\
  \label{eqs24}\widehat{Z^N_{t_k}}=&\frac{1}{h} \E_{t_k}[\Delta W_k
  \widehat{Y^N_{t_{k+1}}}],
\end{align}
where we set $\nabla_x   f^N_{t_k}=\nabla_x f(t_k, X^N_{t_k}, Y^N_{t_{k+1}}, Z^N_{t_k})$ and analogously for $\nabla_y   f^N_{t_k}$ and $\nabla_z   f^N_{t_k}$. Indeed, we can start from (\ref{eqs2}-\ref{eqs3}) and interchange conditional
expectations and weak derivatives (see Proposition 1.2.4 in
\cite{nualart_95}). Another way to get (\ref{eqs23}-\ref{eqs24}) is to take advantage
of the Markov structure of $(X^N_{t_k})_k$ to write $Y^N_{t_k}=y^N(t_k,X^N_{t_k})$,
where the function $y^N$ is the solution of a dynamic programming equation, and then
apply the chain rule. We omit further details.

From (\ref{eqs41}), we also have
\begin{align}\label{eqs42}
  \widehat{Y_t}=\nabla_x u(t,X_t)
  \widehat{X_t},\; \widehat{Z_t}=\nabla_x(\nabla_x
  u\sigma)^*(t,X_t)\widehat{X_t}.
\end{align}
For the sake of clarity, let us write, for any process $V$, $\widehat{\Delta V_k}=D_r
V^N_{t_k}-D_r V_{t_k}$. In particular, we have $\widehat{\Delta
  \overline{Z}_k}=D_r({Z^N_{t_k}}^*-\overline{Z}^*_{t_k})=\widehat{Z^N_{t_k}}-\widehat{\overline{Z}_{t_k}},$ where $\widehat{\overline{Z}_{t_k}}$ is defined as
$h\widehat{\overline{Z}_{t_k}}=\E_{t_k}
\int_{t_k}^{t_{k+1}}\widehat{Z_s}ds$ (see the beginning of Section \ref{section:proof:theo1}).

\subsubsection{Preparatory estimates}
In this part we give some $L_p$-estimates ($p\ge 1$), which are repeatedly used inthe following calculations.
\begin{align}
  \label{eqs33}\bullet\ &\sup_{i \leq j \leq N}(\E_{t_i} |\widehat{X^N_{t_j}}|^{2p})\leq
  C|\widehat{X^N_i}|^{2p},\\
  \label{eqs32}\bullet\ &\E(\max_{0\leq j \leq N} |\widehat{X^N_{t_j}}|^{2p})=O(1),\\ \label{eqs34} \bullet\ &\forall j \in 0..N-1, \; |\widehat{Y^N_{t_j}}|^2\leq
  C|\widehat{X^N_{t_j}}|^2, \;\; \E(\max_{0\leq j \leq N}
  |\widehat{Y^N_{t_j}}|^{2p})=O(1),\\
  \label{eqs45}\bullet\ &\E(\sup_{0 \leq t \leq T}
  |\widehat{X_t}|^{2p}+\sup_{0\leq t \leq T}|\widehat{Y_t}|^{2p}+\sup_{0\leq t \leq T}|\widehat{Z_t}|^{2p})=O(1),\\
  \label{eqs37}\bullet\ &\mbox{Let } F \mbox{ satisfy
  }R_3.\mbox{ Then, } |\E(F(\widehat{X^N_t} -
  \widehat{X_t}))|=O(h).
  \mbox{ Furthermore, } \notag\\
  &\sup_{0\leq k \leq N}\E|\widehat{\Delta X_k}|^{2p}=O(h^p).\\
\bullet\ &\mbox{Analogously to (\ref{eqs18}), $\forall s \in[t_k,t_{k+1}]$, we have}\notag\\
\label{eqs39}
 & \E_{t_k}\big(
  |\widehat{X_s}-\widehat{X_{t_k}}|^{2p}+|\widehat{Y_s}-\widehat{Y_{t_k}}|^{2p}+|\widehat{Z_s}-\widehat{\overline{Z}_{t_k}}|^{2p}\big) =O_k(h^p).
\end{align}
Note that $\widehat{X^N_{t_1}}=\sigma(0,x)$, and $\widehat{X^N_{t_{k+1}}}=(1+h   b'_x(t_k, X^N_{t_k})+\sum_{i=1}^q (\sigma_i)'_x(t_k,X^N_{t_k})\Delta
  W^i_k)\widehat{X^N_{t_k}}$ for $1\le k \le N$. Thus, we easily get $\E_{t_i}|\widehat{X^N_{t_j}}|^{2p}\leq
(1+Ch)\E_{t_i}|\widehat{X^N_{t_{j-1}}}|^{2p}$, and (\ref{eqs33}) follows.
 The proof of (\ref{eqs32}) can be done as the proof of \eqref{eq:b}.

{\it Proof of (\ref{eqs34}).}
From (\ref{eqs23}), we use Young's inequality and boundedness of $\nabla f$ to get
\begin{align}\label{eqs40}
  |\widehat{Y^N_{t_i}}|^2 \leq (1+ \gamma h)|\E_{t_i}
  \widehat{Y^N_{t_{i+1}}}|^2+Ch(h+\frac{1}{\gamma})\big(| \widehat{X^N_{t_i}}|^2+\E_{t_i}  |\widehat{Y^N_{t_{i+1}}}|^2+| \widehat{Z^N_{t_i}}|^2\big).
\end{align}
From (\ref{eqs24}) and the Cauchy Schwarz inequality, we obtain $h|\widehat{Z^N_{t_i}}|^2
\leq C(\E_{t_i}| \widehat{Y^N_{t_{i+1}}}|^2-|\E_{t_i}
\widehat{Y^N_{t_{i+1}}}|^2)$. Hence, with an appropriate choice of $\gamma$, (\ref{eqs40}) is reduced to $ |\widehat{Y^N_{t_i}}|^2 \leq (1+Ch)\E_{t_i}|
  \widehat{Y^N_{t_{i+1}}}|^2+Ch|\widehat{X^N_{t_i}}|^2$, and thus Gronwall's lemma
  yields
  \begin{align*}
    |\widehat{Y^N_{t_i}}|^2 \leq
    C\E_{t_i}(|\widehat{Y^N_{t_N}}|^2+h\sum_{j=i}^{N-1}|\widehat{X^N_{t_j}}|^2) \leq
    C\sup_{i \leq j \leq N-1}  \E_{t_i}|\widehat{X^N_{t_j}}|^2.
  \end{align*}
Finally, estimates (\ref{eqs33}) and (\ref{eqs32}) complete the proof.
  
  {\it Proof of (\ref{eqs45})}. $\E(\sup_{0\leq t \leq T}
  |\widehat{X_t}|^{2p})=O(1)$ follows from (\ref{eqs48}). The other estimates come from this result and (\ref{eqs42}).

  {\it Proof of (\ref{eqs37})}. Let us introduce $X_t'=\nabla_x X_{t} (\nabla_x X_r)^{-1} \sigma(0,x)$ and write $\widehat{X^N_t}-\widehat{X_t}=\widehat{X^N_t}-X_t'+X_t'-\widehat{X_t}$. 

Since $\widehat{X_t}=\nabla_x X_{t} (\nabla_x X_r)^{-1}\sigma(r,X_r)$, a direct application of Proposition \ref{prop5} with $U_t=\sigma(t,X_t)$ gives $\E(F(X_t'-\widehat{X_t}))=O(h)$ for $F$ satisfying $R_2$. Moreover, simple increment estimates yield $\sup_{t\le T}\E|X_t'-\widehat{X_t}|^{2p}=O(h^p)$. 

It remains to study the impact of the difference $\widehat{X^N_t}-X_t'$. $(\widehat{X^N_t})_{t\ge r}$ and $(X_t')_{t\ge r}$ are solutions of 
\begin{align}\label{eqs43}
  \widehat{X^N_t}=&\sigma(0,x)+\int_r^t
  b'_x(\eta(s),X^N_{\eta(s)})\widehat{X^N_{\eta(s)}} ds + \sum_{i=1}^q\int_r^t
  (\sigma_i)'_x(\eta(s),X^N_{\eta(s)})\widehat{X^N_{\eta(s)}} dW^i_s,\notag\\
  X'_t=&\sigma(0,x)+\int_r^t
  b'_x(s,X_s) X'_s ds + \sum_{i=1}^q\int_r^t
  (\sigma_i)'_x(s,X_s) X'_s dW^i_s.
\end{align}
For the sake of simplicity, we take $b\equiv 0$ and $d=q=1$. If we set $\sigma'(s)=\int_0^1\sigma'_x(s, X_s +\lambda(X^N_s -X_s))d\lambda$, we observe that $\Delta X_t$ solves the linear equation $\Delta X_t=\int_0^t [\sigma(\eta(s),X^N_{\eta(s)})-\sigma(s,X^N_s)] dW_s +\int_0^t \sigma'(s) \Delta X_s dW_s,$ which solution is given by (see Theorem $56$ p. $271$ in \cite{protter_90})
\begin{align*}
  \Delta X_t=&\epsilon_t \int_0^t \epsilon_s^{-1}[\sigma(\eta(s),X^N_{\eta(s)})-\sigma(s,X^N_s)]  (dW_s -\sigma'(s) ds)\\
=&-\epsilon_t \int_0^t \epsilon_s^{-1} \big[\int_{\eta(s)}^s \sigma'_x(v,X^N_v)\sigma(\eta(v),X^N_{\eta(v)})dW_v\\
&+ (\sigma'_t(v,X^N_v)+\frac 12\sigma''_{xx}(v,X^N_v)\sigma^2(\eta(v),X^N_{\eta(v)}))dv\big](dW_s -\sigma'(s) ds)
\end{align*} 
where $\epsilon_t= 1+\int_0^t \sigma'(s) \epsilon_s dW_s$. This proves that $\Delta X_t$ satisfies the property \ref{eq:se}. Analogously, if we define $\sigma''(s)=\int_0^1\sigma''_{xx}(s, X_s +\lambda(X^N_s -X_s))d\lambda$ and $\epsilon^N_t= 1+\int_r^t \sigma'_x(s,X^N_s)\epsilon^N_s dW_s$, simple computations lead to
\begin{align*}
 \widehat{X^N_t}-X'_t=&\epsilon^N_t \int_r^t (\epsilon^N_s)^{-1}([\sigma'_x(\eta(s),X^N_{\eta(s)})\widehat{X^N_{\eta(s)}}-\sigma'_x(s,X^N_s)\widehat{X^N_s}]+ \sigma''(s) X'_s \Delta X_s)\\
&\qquad (dW_s -\sigma'_x(s,X^N_s) ds).
\end{align*}
From the above representation, it is straightforward to conclude $\sup_{t\le T}\E|\widehat{\Delta X_t}|^{2p}=O(h^p)$. Now, let us upper bound $\E(F(\widehat{X^N_t}-X'_t))$ which can be decomposed into several terms.
\begin{itemize}
\item The contribution associated to $\epsilon^N_t \int_r^t (\epsilon^N_s)^{-1} [\sigma'_x(\eta(s),X^N_{\eta(s)})\widehat{X^N_{\eta(s)}}-\sigma'_x(s,X^N_s)\widehat{X^N_s}](dW_s -\sigma'_x(s,X^N_s) ds)$ satisfies property \ref{eq:se}, thus Proposition \ref{prop6} yields the expected result.
\item The contribution $\E(F\epsilon^N_t \int_r^t (\epsilon^N_s)^{-1}\sigma''(s) X'_s \Delta X_s \sigma'_x(s,X^N_s) ds)$ is equal to \break $\int_r^t \E(F\epsilon^N_t(\epsilon^N_s)^{-1}\sigma''(s) X'_s \Delta X_s \sigma'_x(s,X^N_s)) ds=O(h)$ in view of Remark  \ref{remarque4}.
\item In the same way, the duality relationship ensures that the last contribution $\E(F\epsilon^N_t \int_r^t (\epsilon^N_s)^{-1}\sigma''(s) X'_s \Delta X_s dW_s)=\int_r^t \E(D_s (F\epsilon^N_t) (\epsilon^N_s)^{-1}\sigma''(s) X'_s \Delta X_s)ds$ is a $O(h)$ (using here that $F$ satisfies $R_3$).
\end{itemize}  

{\it Proof of (\ref{eqs39})}. In view of $\widehat{X_t}=D_r X_{t}= \nabla_x X_{t} (\nabla_x X_r)^{-1}
\sigma(r,X_r)$, the estimate on the increments of $\widehat{X_t}$ becomes clear. The
other ones easily follow.\qed

\subsubsection{Proof of $\max_{1 \leq i \leq N}\E|\widehat{\Delta
    Y_{i}}|^2=O(h).$}\label{section1bis}
Assume that for some non negative random variable $\Lambda_k=O_k(h)+|\Delta X_k|^2 +|\Delta Z_k|^2$, one has
\begin{align}\label{eqs44}
  |\widehat{\Delta Y_k}|^2\leq(1+Ch)\mathbb{E}_{t_k}|\widehat{\Delta
    Y_{k+1}}|^2
  +h |\widehat{\Delta X_k}|^2
  +h \Lambda_kO_k(1).
\end{align}
Take the expectation on both sides, use estimates (\ref{eqs37}) and those of Proposition \ref{theo3} to get 
\begin{align*}
  \E|\widehat{\Delta Y_k}|^2\leq C\E|\widehat{\Delta Y_N}|^2+ O(h)+C h\sum_{k=0}^{N-1}\E(|\Delta
  Z_k|^2 O_k(1)).
\end{align*}
On the one hand, as $\widehat{\Delta Y_N}=\Phi'(X^N_{t_N})\widehat{X^N_{t_N}}-\Phi'(X_{t_N})\widehat{X_{t_N}}$, clearly $\E|\widehat{\Delta Y_N}|^2=O(h)$. On the other hand, in view of \eqref{eqs10} with $p=2$, the summation above is a $O(h)$. This proves $\max_{1 \leq k \leq N}\E|\widehat{\Delta Y_{k}}|^2=O(h).$

{\bf Proof of (\ref{eqs44}).}
From (\ref{eqs22}) and (\ref{eqs23}), we obtain
\begin{align*}
  \widehat{\Delta Y_k}=\E_{t_k}(\widehat{\Delta Y_{k+1}})+\E_{t_k}&\big(\int_{t_k}^{t_{k+1}} [\nabla_x
  f^N_{t_k}\widehat{X^N_{t_k}}-
  f'_x(\theta_s)\widehat{X_s}\\
  &+\nabla_y
  f^N_{t_k}\widehat{Y^N_{t_{k+1}}}- f'_y(\theta_s)\widehat{Y_s}
  +\nabla_z f^N_{t_k}\widehat{Z^N_{t_k}}-
  f'_z(\theta_s)\widehat{Z_s}]ds\big).
\end{align*}
Since $f\in C^{2,4,4,4}_b$, it follows that for any $\gamma>0$ (to be fixed later)
\begin{align}\label{eqs60}
  |\widehat{\Delta Y_k}|^2\leq& (1+\gamma h)|\E_{t_k}(\widehat{\Delta
    Y_{k+1}})|^2+C(h+\frac{1}{\gamma })\E_{t_k}\big(\int_{t_k}^{t_{k+1}} [|\nabla_x
  f^N_{t_k}\widehat{X^N_{t_k}}-
  f'_x(\theta_s)\widehat{X_s}|^2\notag\\
  &+|\nabla_y
  f^N_{t_k}\widehat{Y^N_{t_{k+1}}}- f'_y(\theta_s)\widehat{Y_s}|^2
  +|\nabla_z f^N_{t_k}\widehat{Z^N_{t_k}}-f'_z(\theta_s)\widehat{Z_s}|^2] ds \big)\\
\label{eqs25}
\leq& (1+\gamma h)|\E_{t_k}(\widehat{\Delta
    Y_{k+1}})|^2+C(h+\frac{1}{\gamma}) (T_k^1+T_k^2),
\end{align}
where we put $T_k^1=\E_{t_k}(\int_{t_k}^{t_{k+1}} [|\widehat{X^N_{t_k}}-\widehat{X_s}|^2  +|\widehat{Y^N_{t_{k+1}}}-\widehat{Y_s}|^2
  +|\widehat{Z^N_{t_k}}-\widehat{Z_s}|^2] ds )$, $T_k^2=\E_{t_k}\int_{t_k}^{t_{k+1}}(h+|X_s-X^N_{t_k}|^2+|Y_s-Y^N_{t_{k+1}}|^2+|Z_s-Z^N_{t_k}|^2)(|\widehat{X_s}|^2+|\widehat{Y_s}|^2+|\widehat{Z_s}|^2) ds$. To get (\ref{eqs44}), we need to simplify (\ref{eqs25}), by estimating $T_k^1$ and $T_k^2$.\\
{\bf Term $T_k^1$}. Firstly, we write $\E_{t_k} |\widehat{Y^N_{t_{k+1}}}-\widehat{Y_s}|^2 \leq 2\E_{t_k}|\widehat{Y_{t_{k+1}}}-\widehat{Y_s}|^2
+2\E_{t_k}|\widehat{\Delta Y_{k+1}}|^2$. We do the same for
$\widehat{X^N_{t_k}}-\widehat{X_s}$. Then, the usual increment estimates yield
\begin{align*}
  \E_{t_k}
|\widehat{Y^N_{t_{k+1}}}-\widehat{Y_s}|^2+\E_{t_k}|\widehat{X^N_{t_k}}-\widehat{X_s}|^2
\leq O_k(h)+2|\widehat{\Delta X_k}|^2 +2\E_{t_k}|\widehat{\Delta Y_{k+1}}|^2.
\end{align*}
 Secondly, analogously to (\ref{eqs6}), we have
\begin{align*}
  \mathbb{E}_{t_k}   \int_{t_k}^{t_{k+1}} |\widehat{Z^N_{t_k}}-\widehat{Z_s}|^2 ds = \mathbb{E}_{t_k}
  \int_{t_k}^{t_{k+1}} |\widehat{\overline{Z}_{t_k}}-\widehat{Z_s}|^2 ds + h \mathbb{E}_{t_k}|\widehat{Z_{t_k}^N}-\widehat{\overline{Z}_{t_k}}|^2. 
\end{align*}
Finally, we obtain $T_k^1\leq 
  C h( O_k(h)+|\widehat{\Delta X_k}|^2
  +\E_{t_k}|\widehat{\Delta Y_{k+1}}|^2+|\widehat{\Delta\overline{Z}_k}|^2).$

{\bf Term $T_k^2$.} Easy calculations combining \eqref{eqs18}, Proposition \ref{prop2} and \eqref{eqs45} give $ T_k^2 \leq(O_k(h^2)+h|\Delta X_k|^2+h|\Delta Z_k|^2)O_k(1)= h \Lambda_k O_k(1).$

{\bf Conclusion.} Plugging the estimates on $T^1_k$ and $T^2_k$ into (\ref{eqs25}), we get
\begin{align}\label{eqs46}
  |\widehat{\Delta Y_k}|^2\leq& (1+\gamma h)|\E_{t_k}(\widehat{\Delta
    Y_{k+1}})|^2+Ch(h+\frac{1}{\gamma})|\widehat{\Delta \overline{Z}_k}|^2\notag\\
  &+Ch(h+\frac{1}{\gamma})(|\widehat{\Delta X_k}|^2
  +\E_{t_k}|\widehat{\Delta Y_{k+1}}|^2+\Lambda_kO_k(1)).
\end{align}
Note that $h\widehat{\overline{Z}_{t_k}}=\E_{t_k}(\Delta W_k(\widehat{Y_{t_{k+1}}}+\int_{t_k}^{t_{k+1}}  [f'_x(\theta_s)\widehat{X_s}  + f'_y(\theta_s)\widehat{Y_s} +f'_z(\theta_s)\widehat{Z_s}] ds))$, whence $h\widehat{\Delta \overline{Z}_{t_k}}=\E_{t_k}(\Delta W_k(\widehat{\Delta
  Y_{k+1}}+\int_{t_k}^{t_{k+1}} [ f'_x(\theta_s)
\widehat{X_s}  + f'_y(\theta_s)\widehat{Y_s} +f'_z(\theta_s)\widehat{Z_s}]
ds))$. By proceeding as before, we easily prove 
\begin{align}\label{eqs27}
  h|\widehat{\Delta \overline{Z}_{t_k}}|^2 \leq C(
  \mathbb{E}_{t_k}|\widehat{\Delta Y_{k+1}}|^2-
  |\mathbb{E}_{t_k}\widehat{\Delta Y_{k+1}}|^2)+O_k(h^2).
\end{align}
Combining this upper bound with \eqref{eqs46} for a good choice of $\gamma$ gives (\ref{eqs44}).\qed
\subsubsection{Proof of $h\E\big(\sum_{k=1}^{N-1}|\widehat{\Delta Z_k}|^2
  \big)=O(h).$}
In view of (\ref{eqs39}), this is equivalent to prove $h\E\big(\sum_{k=1}^{N-1}|\widehat{\Delta \overline{Z}_k}|^2
  \big)=O(h).$ To establish this estimate, we start from (\ref{eqs27}) to get
\begin{align}\label{eqs29}
  h \sum_{k=1}^{N-1}\E|\widehat{\Delta \overline{Z}_k}|^2 \leq C\sum_{k=1}^{N-1}(
  \E|\widehat{\Delta Y_k}|^2-\E
  |\mathbb{E}_{t_k}\widehat{\Delta Y_{k+1}}|^2)+C\E|\widehat{\Delta Y_N}|^2+O(h).
\end{align}
Now, we work on $|\widehat{\Delta Y_k}|^2-|\mathbb{E}_{t_k}\widehat{\Delta
  Y_{k+1}}|^2$.
The choice $\gamma=2C^2$ in (\ref{eqs46}) leads to
\begin{align*}
  |\widehat{\Delta Y_k}|^2-&|\E_{t_k}(\widehat{\Delta
    Y_{k+1}})|^2\leq \gamma h|\E_{t_k}(\widehat{\Delta
    Y_{k+1}})|^2+h(\frac{1}{2C}+Ch)|\widehat{\Delta \overline{Z}_k}|^2\\
  &+h(Ch+\frac{1}{2C})(|\widehat{\Delta X_k}|^2
  +\E_{t_k}|\widehat{\Delta Y_{k+1}}|^2+\Lambda_kO_k(1)).
\end{align*}
From (\ref{eqs37}) and the result from Section \ref{section1bis}, we have $\max_{1 \le k \le N}\E(|\widehat{\Delta X_k}|^2
  +|\widehat{\Delta Y_{k}}|^2)=O(h)$. We also have
  $\E(\Lambda_kO_k(1))=O(h)+\E(|\Delta Z_k|^2 O_k(1))$. Consequently, for $h$ small enough, one has $\E|\widehat{\Delta Y_k}|^2-\E|\E_{t_k}(\widehat{\Delta
    Y_{k+1}})|^2\leq  \frac{2h}{3C}\E|\widehat{\Delta
    \overline{Z}_k}|^2+O(h^2)+ Ch\E(|\Delta Z_k|^2 O_k(1)).$
Putting this estimate into (\ref{eqs29}) yields
\begin{align*}
  \frac{1}{3}h\sum_{k=1}^{N-1} \E|\widehat{\Delta \overline{Z}_k}|^2 \leq O(h)+Ch\sum_{k=1}^{N-1}\E(|\Delta Z_k|^2 O_k(1)).
\end{align*}
Inequality \eqref{eqs10} with $p=2$ directly shows that the sum above is a $O(h)$. \qed
\subsection{Expansion of $Z^N_{t_k}-Z_{t_k}$}
We recall that $u \in C^{2,4}_b$ owing to Hypothesis \ref{hyp5}. From (\ref{eqs14}), we have $\Delta Z_k =
O(h)+\frac{1}{h}\E_{t_k}[(u^N(t_{k+1},X^N_{t_{k+1}})-u(t_{k+1},X_{t_{k+1}}))\Delta
W_k^*]$. Let $(X_t^{s,\overline{x}})_{t\ge s}$ denote the solution of the
SDE (\ref{eqs4}) starting at time $s$ from $\overline{x}$. We write $X_t$ for
$X_t^{0,x}.$ Note that $X_{t_{k+1}}=X_{t_{k+1}}^{t_k,X_{t_k}}$. In the same way, the
Euler scheme starting at time $t_k$ at $\overline{x}$ is denoted by
$(X^{N,t_k,\overline{x}}_{t_j})_{j \geq k}$. With this notation we can rewrite $\Delta Z_k$
\begin{align}\label{eqs16}
  \Delta Z_k =&
  \frac{1}{h}\E_{t_k}[(u^N(t_{k+1},X^{N,t_k,X_{t_k}^N}_{t_{k+1}})-u(t_{k+1},X_{t_{k+1}}))\Delta
  W_k^*]+O(h),\notag\\
  =&\frac{1}{h}\E_{t_k}[(u(t_{k+1},X^{t_k,X_{t_k}^N}_{t_{k+1}})-u(t_{k+1},X_{t_{k+1}}))\Delta W_k^*]\notag\\
  &+\frac{1}{h}\E_{t_k}[(u^N(t_{k+1},X^{N,t_k,X_{t_k}^N}_{t_{k+1}})-u(t_{k+1},X^{t_k,X_{t_k}^N}_{t_{k+1}}))\Delta W_k^*] +O(h).
\end{align}
We work on the first two terms separately by proving
\begin{lem}\label{lemme4}
    $\frac{1}{h}\E_{t_k}[(u(t_{k+1},X^{t_k,X_{t_k}^N}_{t_{k+1}})-u(t_{k+1},X_{t_{k+1}}))\Delta W_k^*]=O(|\Delta X_k|^2)+O(h)\break +[\nabla_x(\nabla_x
     u\ \sigma)^*(t_k,X_{t_k})\Delta X_k]^*.$
\end{lem}
\begin{lem}\label{lemme5}
   $ \frac{1}{h}\big|\E_{t_k}[(u^N(t_{k+1},X^{N,t_k,X^N_{t_k}}_{t_{k+1}})-u(t_{k+1},X^{t_k,X^N_{t_k}}_{t_{k+1}}))\Delta W_k^*]\big| =O_k(h).$
\end{lem}
The combination of these Lemmas completes the proof of Theorem \ref{theo4}.
\subsubsection{Proof of Lemma \ref{lemme4}.}
For the sake of simplicity, let $\Delta_N
X_{k+1}$ denote $X^{t_k,X_{t_k}^N}_{t_{k+1}}-X_{t_{k+1}}$ (which is different from
$\Delta X_{k+1}=X^{N,t_k,X^N_{t_k}}_{t_{k+1}}-X_{t_{k+1}}$). From a Taylor-Lagrange formula, we obtain 
\begin{align*}
  & u(t_{k+1},X^{t_k,X_{t_k}^N}_{t_{k+1}})-u(t_{k+1},X_{t_{k+1}})=
  u'_x(t_{k+1},X_{t_{k+1}})\Delta_N X_{k+1}\\ 
  &+\int_0^1(1-\lambda)
  (\Delta_N X_{k+1})^*H_x(u)\big(t_{k+1},X_{t_{k+1}}+\lambda\Delta_N
  X_{k+1}\big)\Delta_N X_{k+1} d\lambda. 
\end{align*}
Thus, using the duality relationship, one has
\begin{align*}
  &\E_{t_k}[(u(t_{k+1},X^{t_k,X_{t_k}^N}_{t_{k+1}})-u(t_{k+1},X_{t_{k+1}}))\Delta W_k^*]\nonumber\\&=\int_{t_k}^{t_{k+1}} R^1_k(t) \mbox{dt}+\int_{t_k}^{t_{k+1}} R^2_k(t) \mbox{dt}+ \int_0^1(1-\lambda) R^3_k(\lambda)  d\lambda ,\\
\mbox{with }&R^1_k(t)=\E_{t_k}[(\Delta_N X_{k+1})^* 
H_x(u)(t_{k+1},X_{t_{k+1}})D_tX_{t_{k+1}}],\\
&R^2_k(t)=\E_{t_k}[u'_x(t_{k+1},X_{t_{k+1}})D_t(\Delta_N X_{k+1})],\\
&R^3_k(\lambda)=\E_{t_k}[
  (\Delta_N X_{k+1})^*H_x(u)(t_{k+1},X_{t_{k+1}}+\lambda\Delta_N
  X_{k+1})\Delta_N X_{k+1}\Delta W_k^*].
\end{align*}
{\bf Expansion of $R^1_k(t)$.}
Clearly $\Delta_NX_{k+1}=\Delta X_k+U_{t_{k+1}}-U_{t_k}$, where $U$ is an It\^o process with drift term $\alpha_s=b(s,X^{t_k,X_{t_k}^N}_s)-b(s,X_s)$ and diffusion term $\beta_s=\sigma(s,X^{t_k,X_{t_k}^N}_s)-\sigma(s,X_s)$, both being bounded. Thus, we can apply Proposition \ref{prop5}, letting $F=H_x(u)(t_{k+1},X_{t_{k+1}}) D_tX_{t_{k+1}}$. Because $u\in C^{2,4}_b$ and in view of (\ref{eqs48}), we get
 \[ R^1_k(t)=O(h)+(\Delta  X_k)^*\E_{t_k}[H_x(u)(t_{k+1},X_{t_{k+1}})D_tX_{t_{k+1}}].\]
We expand the latter factor. As $D_t X_{t_{k+1}}= \nabla_x X_{t_{k+1}} (\nabla_x X_t)^{-1} \sigma(t,X_t)$, we have
\begin{align*}
  H_x(u)(t_{k+1},X_{t_{k+1}})D_tX_{t_{k+1}}
  &=(H_x(u)(t_{k+1},X_{t_{k+1}})\sigma(t,X_t)-
  H_x(u)(t,X_t)\sigma(t,X_t))\notag\\
  &+(H_x(u)(t,X_t)\sigma(t,X_t)-
  H_x(u)(t_k,X_{t_k})\sigma(t_k,X_{t_k}))\notag\\
  &+(H_x(u)(t_{k+1},X_{t_{k+1}})[\nabla_x X_{t_{k+1}} (\nabla_x
  X_t)^{-1}-I]\sigma(t,X_t))\notag\\
  &+H_x(u)(t_k,X_{t_k})\sigma(t_k,X_{t_k}).
\end{align*}
The first three contributions in the r.h.s. above can be handled in the same way and we give a detailed proof only for the first one. It is enough to apply Proposition \ref{prop5} with $F=\sigma(t,X_t)$ and $U_s=H_x(u)(s,X_s)$. Then, $\E_{t_k}[F(U_{t_{k+1}}-U_t)]$ is of order $h$ with a constant involving $b, \sigma, u$ and its derivatives up to order $4$. Finally, this gives 
\[ R^1_k(t)=O(h)+(\Delta  X_k)^*H_x(u)(t_k,X_{t_k})\sigma(t_k,X_{t_k}),\]
uniformly in $t\in[t_k,t_{k+1}]$.

{\bf Expansion of $R^2_k(t)$.} For $t_k \le t\leq t_{k+1}$, we have
\begin{align*}
  D_t(\Delta_N X_{k+1})=&[\nabla_x X_{t_{k+1}}^{X_{t_k}^N,t_k} (\nabla_x
  X_t^{X_{t_k}^N,t_k})^{-1}-I]\sigma(t,X_t^{X_{t_k}^N,t_k})\notag\\
  &-[\nabla_x X_{t_{k+1}} (\nabla_x
  X_t)^{-1}-I]\sigma(t,X_t)-(\sigma(t,X_t)-\sigma(t_k,X_{t_k}))\notag\\
  &+\sigma(t,X_t^{X_{t_k}^N,t_k})-\sigma(t_k,X_{t_k}^N)+\sigma(t_k,X_{t_k}^N)-\sigma(t_k,X_{t_k}).
\end{align*}
As before, apply Proposition \ref{prop5} to each of these terms but the last one, with $F=u'_x(t_{k+1},X_{t_{k+1}})$, using $u, b, \sigma \in C^{2,4}_b$ and (\ref{eqs48}). It follows that $R^2_k(t)=O(h)+\E_{t_k}[u'_x(t_{k+1},X_{t_{k+1}})](\sigma(t_k,X_{t_k}^N)-\sigma(t_k,X_{t_k})).$ An application of Itô's formula yields 
\begin{align*}
 R^2_k(t)&=O(h)
  +\sum_{i=1}^d  u'_{x_i}(t_k,X_{t_k})(\sigma^i(t_k,X_{t_k}^N)-\sigma^i(t_k,X_{t_k}))\\
&=O(h+|\Delta X_k|^2)
  +\sum_{i=1}^d  u'_{x_i}\nabla_x
  ([\sigma^i]^*)(t_k,X_{t_k})\Delta X_k,
\end{align*}
uniformly in $t\in[t_k,t_{k+1}]$. Finally, simple matrix computations lead to 
\[R^1_k(t)+R^2_k(t)=O(h+|\Delta X_k|^2)+[\nabla_x(\nabla_x
  u\sigma)^*(t_k,X_{t_k})\Delta X_k]^*.\]

{\bf Upper bound for $R^3_k(\lambda)$.} To complete the proof of Lemma \ref{lemme4}, note that it remains to justify that $R^3_k(\lambda)=hO(h+|\Delta X_k|^2)$ uniformly in $\lambda$. The duality formula gives
\begin{align*}
R^3_k(\lambda)=  &\E_{t_k}[\int_{t_k}^{t_{k+1}}D_t[(\Delta_N X_{k+1})^*H_x(u)(t_{k+1},X_{t_{k+1}}+\lambda\Delta_N
  X_{k+1})\Delta_N X_{k+1}]\mbox{dt}.
\end{align*}
The term in the integral equals $\sum_{i,j=1}^d [2 D_t (\Delta_N
X_{k+1,i})\Delta_N X_{k+1,j} \partial_{x_i,x_j}^2 u(t_{k+1},X_{t_{k+1}}+\lambda\Delta_N  X_{k+1})+ \Delta_N X_{k+1,i} \Delta_N X_{k+1,j}
D_t(\partial_{x_i,x_j}^2 u(t_{k+1},X_{t_{k+1}}+\lambda\Delta_N  X_{k+1}))]$. Thanks to (\ref{eqs48}) and (\ref{eqs49}) and successive applications of Proposition \ref{prop5}, we finally prove our assertion. We omit further details.\qed

\subsubsection{Proof of Lemma \ref{lemme5}}
As for Lemma \ref{lemme1}, we only do the proof for $t_k=0$, i.e. we have to show $ |\E_{t_k}[(u^N(t_1,X^{N,0,x}_{t_1})-u(t_1,X^{0,x}_{t_1}))\Delta W_0^*]| \leq K(T,x)h^2.$ We have
$\E[(u^N(t_1,X^{N,0,x}_{t_1})-u(t_1,X^{0,x}_{t_1}))\Delta
W_0^*]=\E[\Delta Y_1 \Delta W_0^*].$
By using (\ref{eq:y:dl}), we come up with
\begin{align*}
  \E[\Delta Y_1 \Delta W_0^*]=\E[\xi_1...\xi_{N-1}\Delta Y_N
  \Delta W_0^*]+\E[h\sum_{i=1}^{N-1}(f'_x(\theta_{t_i})\Delta X_i+ \tilde{\chi_i}) \xi_1...\xi_{i-1}\Delta W_0^*],
\end{align*}
where $\tilde{\chi_i}=\E_{t_i}(\chi_i)$ ($\xi_i$ and $\chi_i$ are defined in \eqref{eq:xi} and \eqref{eq:chi}).
In the following $\tilde{\eta_i}$ denotes $\xi_1...\xi_{i-1}$ and
$\tilde{\eta_1}=1$. We easily prove that $(\tilde{\eta}_i)_{1 \leq i \leq N}$ has the analogous properties to $(\eta_i)_{0 \leq i \leq N}$. Estimates (\ref{eq:a}) and (\ref{eq:b}) remain valid for $\tilde\eta$ and under Hypothesis \ref{hyp5}, the estimate (\ref{eq:a}) becomes
\begin{equation}
  \label{eq:abis}\tilde{\eta_k} \text{ satisfies } R_3 \text{ uniformly in } k.
\end{equation}

{\bf Step 1 : Proof of $\E[\xi_1...\xi_{N-1}\Delta Y_N
  \Delta W_0^*]=\E[\tilde{\eta_N}\Delta Y_N\Delta W_0^*]=O(h^2)$.}\\
As before, we use the duality formula:
\begin{align*}
  \E[\tilde{\eta_N}\Delta Y_N
  \Delta W_0^*]=&\E \int_{0}^{t_1} (D_t[\tilde{\eta_N}]\Delta
  Y_N+\tilde{\eta_N}D_t[\Delta Y_N])\mbox{dt}.
\end{align*}
Since $ \tilde{\eta_N}$ satisfies (\ref{eq:abis}), we proceed as in {\bf Step 4} of Lemma
\ref{lemme1} and we get $\E(D_t[\tilde{\eta_N}]\Delta
Y_N)=O(h).$
Furthermore, we have
\[  D_t[\Delta Y_N]=( \Phi'(X^N_T)-\Phi'(X_T))D_t
  X^N_T+ \Phi'(X_T)(D_t X^N_T-D_t X_T).\]
On the one hand, analogously to previous computations, we establish  $\E(\tilde{\eta_N}(\Phi'(X^N_T)-\Phi'(X_T))D_t
X^N_T)=O(h)$.\\
 On the other hand, we prove $\E(\tilde{\eta_N}\Phi'(X_T)(D_t X^N_T-D_t X_T))=O(h)$. Thanks to (\ref{eqs48}) and
 (\ref{eq:a}), $\tilde{\eta_N}\Phi'(X_T)$ satisfies condition $R_3$. Then, by
 applying (\ref{eqs37}), we get the result.
 
{\bf Step 2 : Proof of $\E[h\sum_{i=1}^{N-1}f'_x(\theta_{t_i})\Delta X_i
  \xi_1...\xi_{i-1}\Delta W_0^*]=O(h^2)$.}\\
This is a similar proof to the one done at {\bf Step 1}, with $\Phi(x)=x$.

{\bf Step 3 : Proof of $\E [h\sum_{i=1}^{N-1} \tilde{\chi_i}
  \tilde{\eta_i}\Delta W_0^*]=O(h^2).$}\\
A careful inspection of the definition of $G_0, G_y$ and $G_z$ appearing in \eqref{eq:chi} shows that under Hypothesis \ref{hyp5}, these functions are continuously differentiable w.r.t. the variable $x$ (with a bounded derivative). Hence, if we write $\chi_i=\chi^1_i+\int_0^1 (1-\lambda)\chi^2_i(\lambda)d\lambda$ with (see \eqref{eq:chi})
\begin{align*}
\chi^1_i&=\int_{t_i}^{t_{i+1}} (G_0(s,X_s)+  f_y'(\theta_{t_i})G_y(s,X_s)+ f_z'(\theta_{t_i})G_z(s,X_s))ds,\\
\chi^2_i(\lambda)&= \Delta X_i^* f^{''}_{xx}(\theta^{\lambda}_{t_i})\Delta    X_i\break +f^{''}_{yy}(\theta^{\lambda}_{t_i})(Y^N_{t_{i+1}}-Y_{t_i})^2+\Delta Z_i    f^{''}_{zz}(\theta^{\lambda}_{t_i})\Delta Z_i^*    \\&+2\Delta X_i^*f^{''}_{xy}(\theta^{\lambda}_{t_i})(Y^N_{t_{i+1}}-Y_{t_i})\break +2\Delta X_i^*f^{''}_{xz}(\theta^{\lambda}_{t_i})\Delta Z_i^*+ 2(Y^N_{t_{i+1}}-Y_{t_i})f^{''}_{yz}(\theta^{\lambda}_{t_i})\Delta Z_i^*,
\end{align*}
we note that the random variable $\chi_i$ is in $\mathbb{D}^{1,\infty}$. Thus and because $\tilde{\chi}_i=\E_{t_i}(\chi_i)$, one has $\E[\tilde{\chi_i} \tilde{\eta_i}\Delta W_0^*]=\E[\chi_i \tilde{\eta_i}\Delta W_0^*]=\E[ \int_0^{t_1} (\chi_i D_t\tilde \eta_i+\tilde\eta_i D_t \chi_i) dt]$.

The upper bound $\tilde{\chi}_i=\E_{t_i}(\chi_i)=O_i(h) + O(|\Delta X_i|^2+ |\Delta
Z_i|^2)$ (see \eqref{eq:chik}) is sufficient to show $\E[\sum_{i=1}^{N-1}\chi_i D_t
\tilde\eta_i]=O(1)$ uniformly in $t$ (follow the arguments of the conclusion of the
proof of Lemma \ref{lemme1} and use (\ref{eq:b}) with $\tilde \eta$).

Now, it remains to establish $\E[\sum_{i=1}^{N-1}\tilde\eta_i D_t\chi_i]=O(1)$. On the one hand, clearly $\E_{t_i}[D_t\chi_i^1]=O_i(h)$ and we conclude 
$\E[\sum_{i=1}^{N-1}\tilde\eta_i D_t\chi^1_i]=O(1)$ uniformly in $t$. On the other
hand, $\chi^2_i$ can be decomposed into several contributions, which can be analyzed
with the same arguments. Let us detail how to handle one of them, for instance
$\E[\sum_{i=1}^{N-1}\tilde{\eta_i}D_t(\Delta
X_i^*f^{''}_{xz}(\theta^{\lambda}_{t_i})\Delta Z_i^*)]$ which has to be a $O(1)$. We do the proof for $d=q=1$. Write $D_t(\Delta
X_if^{''}_{xz}(\theta^{\lambda}_{t_i})\Delta Z_i)=\Delta
X_if^{''}_{xz}(\theta^{\lambda}_{t_i})D_t(\Delta Z_i)+D_t(\Delta
X_i)f^{''}_{xz}(\theta^{\lambda}_{t_i})\Delta Z_i+\Delta
X_iD_t(f^{''}_{xz}(\theta^{\lambda}_{t_i}))\Delta Z_i $. As $f^{''}$ is bounded, we have
\begin{align*}
  \big|\E[\sum_{i=1}^{N-1}\tilde{\eta_i}\Delta
X_if^{''}_{xz}(\theta^{\lambda}_{t_i})D_t(\Delta Z_i)]\big|&\leq\E\big[\sum_{i=1}^{N-1}|\tilde{\eta_i}||\Delta
X_i||f^{''}_{xz}(\theta^{\lambda}_{t_i})||D_t(\Delta Z_i)|\big]\\
&\leq C\big(\E(\sum_{i=1}^{N-1}(|\tilde{\eta_i}|^2|\Delta
X_i|^2))\big)^{\frac{1}{2}}\big(\E(\sum_{i=1}^{N-1}|D_t(\Delta
Z_i)|^2)\big)^{\frac{1}{2}}.
\end{align*}
Thanks to Proposition \ref{lemme3}, (\ref{eq:abis}) and Proposition \ref{theo3}, we get that \linebreak[4] $\E[\sum_{i=1}^{N-1}\tilde{\eta_i}\Delta
X_if^{''}_{xz}(\theta^{\lambda}_{t_i})(D_t \Delta Z_i)]=O(1)$. Analogously, using (\ref{eq:abis}-\ref{eqs10}-\ref{eqs37}), we obtain
$\E[\sum_{i=1}^{N-1}\tilde{\eta_i}(D_t \Delta
X_i)f^{''}_{xz}(\theta^{\lambda}_{t_i})\Delta Z_i]=O(1).$
It remains to demonstrate that $ \big|\E[\sum_{i=1}^{N-1}\tilde{\eta_i}\Delta
X_iD_t(f^{''}_{xz}(\theta^{\lambda}_{t_i}))\Delta Z_i]\big|=O(1).$ We have
\begin{align*}
&   D_t(f^{''}_{xz}(\theta^{\lambda}_{t_i}))=f^{'''}_{xzx}(\theta^{\lambda}_{t_i})(\lambda D_t X^N_{t_i}+(1-\lambda)D_t X_{t_i})\\
 & +f^{'''}_{xzy}(\theta^{\lambda}_{t_i})(\lambda D_t Y^N_{t_{i+1}}+(1-\lambda)D_t Y_{t_i})+f^{'''}_{xzz}(\theta^{\lambda}_{t_i})(\lambda D_t Z^N_{t_i}+(1-\lambda)D_t Z_{t_i}).
\end{align*}
The most difficult term to bound among these three ones is the
one which contains $D_t Z^N_{t_i}$.
If we write $\lambda D_t Z^N_{t_i}+(1-\lambda)D_t Z_{t_i}=\lambda D_t(\Delta
Z_i)+D_t Z_{t_i}$, we obtain
\begin{align*}
  &\big|\E[\sum_{i=1}^{N-1}\tilde{\eta_i}\Delta
  X_if^{'''}_{xzz}(\theta^{\lambda}_{t_i})\lambda D_t(\Delta Z_i) \Delta Z_i]\big|\\
&  \leq C\big(\E(\sum_{i=1}^{N-1}|D_t(\Delta
Z_i) |^2)\big)^{\frac{1}{2}}(\E(\sum_{i=1}^{N-1}(|\Delta X_i|^2|\tilde{\eta_i}|^2|\Delta
Z_i|^2))\big)^{\frac{1}{2}},\\
&\leq C\big(\E(\sum_{i=1}^{N-1}|D_t(\Delta
Z_i) |^2)\big)^{\frac{1}{2}}\big(\E(\sum_{i=1}^{N-1}|\Delta
Z_i|^2)^2\big)^{\frac{1}{4}}\big(\E(\max_{0\leq i \leq N}|\tilde{\eta_i}|^4\max_{0\leq i \leq N}|\Delta
X_i|^4)\big)^{\frac{1}{4}},
\end{align*}
 Applying Proposition \ref{lemme3}, (\ref{eqs10}), Proposition \ref{theo3} and (\ref{eq:b}) (with $\tilde \eta$)
 lead to \linebreak[4]$\E[\sum_{i=1}^{N-1}\tilde{\eta_i}\Delta
 X_if^{'''}_{xz}(\theta^{\lambda}_{t_i})\lambda D_t(\Delta Z_i) \Delta
 Z_i]=O(1).$
 Proposition \ref{theo3}, (\ref{eqs10}), (\ref{eqs32}), (\ref{eqs34}), (\ref{eqs45})
 and (\ref{eq:abis}) enable us to prove that the others terms of \linebreak[4]$\E[\sum_{i=1}^{N-1}\tilde{\eta_i}\Delta
X_iD_t(f^{''}_{xz}(\theta^{\lambda}_{t_i}))\Delta Z_i]$ are $O(1)$.\qed
\bibliographystyle{elsart-harv.bst}
\bibliography{biblio}

\end{document}